\documentclass[12pt,reqno]{amsproc}
\usepackage{graphicx}
\usepackage{pstricks,pst-plot}
\usepackage{psfrag}
\usepackage{amsmath,amssymb}
\usepackage{afterpage}

\newtheorem{theorem}{Theorem}
\newtheorem{lemma}[theorem]{Lemma}

\newtheorem{definition}[theorem]{Definition}

\newtheorem{remark}[theorem]{Remark}
\newtheorem{example}[theorem]{Example}

\newcounter{listacnt}\renewcommand{\thelistacnt}{\alph{listacnt}}

\newcommand{\rmref}[1]{{Eqn.~\rm(\ref{#1})}}

\newcommand{\abs}[1]{{\left\vert{#1}\right\vert}}
\newcommand{\norm}[1]{{\|{#1}\|}}

\newcommand{\re}{{\mathbb{R}}}
\newcommand{\real}{{\mathbb{R}}}
\newcommand{\zet}{{\mathbb{Z}}}
\newcommand{\zed}{{\mathbb{Z}}}
\newcommand{\nat}{{\mathbb{N}}}
\newcommand{\queu}{{\mathbb{Q}}}

\renewcommand{\u}{{\bf u}}
\renewcommand{\v}{{\bf v}}
\newcommand{\w}{{\bf w}}
\newcommand{\z}{{\bf z}}
\newcommand{\ww}{{\bf w}}
\newcommand{\rel}{~{\mbox{\sc rel}}~}
\newcommand{\R}{{\mathcal R}}

\newcommand{\A}{{\mathcal A}}
\newcommand{\uu}{{\bf u}}
\newcommand{\vv}{{\bf v}}

\newcommand{\inv}{{\rm Inv}}

\newcommand{\HH}{{\bf H}}
\newcommand{\Conf}{{\Omega}}
\newcommand{\DConf}{{\mathcal D}}
\newcommand{\pl}{\mbox{\sc pl}}
\renewcommand{\sp}{\mbox{\sc sp}}
\newcommand{\disc}{\mbox{\sc disc}}
\newcommand{\length}{{\iota}}
\renewcommand{\d}{\Delta}
\newcommand{\dd}{\Delta^2}
\newcommand{\EE}{{\mathcal E}}
\newcommand{\del}{\partial}
\newcommand{\cl}{\mbox{\sc cl}}

\newcommand{\dist}{\text{\textup{dist}}}

\newgray{grey1}{0.75}\newgray{grey2}{0.6}\newgray{grey3}{0.45}

\begin{document}

\title[SCALAR PARABOLIC PDE'S AND BRAIDS]
{SCALAR PARABOLIC PDE'S AND BRAIDS}

\author{{\sc R. W. Ghrist}}
\address{Department of Mathematics,
University of Illinois, Urbana, IL 61801, USA}
\email{ghrist@math.uiuc.edu}
\thanks{RG supported in part by NSF CAREER grant DMS-0337713.}

\author{{\sc R. C. Vandervorst}}
\address{Department of Mathematics, Vrije Universiteit Amsterdam,
De Boelelaan 1081, 1081 HV, Amsterdam, Netherlands.}
\email{vdvorst@few.vu.nl}
\thanks{RCV supported in part by NWO VIDI grant 639.032.202 and
RTN grant HPRN-CT-2002-00274.  \today}

\date{\today}
\maketitle

\begin{abstract}
The comparison principle for scalar second order parabolic PDEs on
functions $u(t,x)$ admits a topological interpretation: pairs of
solutions, $u^1(t,\cdot)$ and $u^2(t,\cdot)$, evolve so as to not
increase the intersection number of their graphs. We generalize to
the case of multiple solutions
$\{u^\alpha(t,\cdot)\}_{\alpha=1}^n$. By lifting the graphs to
Legendrian braids, we give a global version of the comparison
principle: the curves $u^\alpha(t,\cdot)$ evolve so as to (weakly)
decrease the algebraic length of the braid.

We define a Morse-type theory on Legendrian braids which we
demonstrate is useful for detecting stationary and periodic
solutions to scalar parabolic PDEs. This is done via
discretization to a finite dimensional system and a suitable
Conley index for discrete braids.

The result is a toolbox of purely topological methods for finding
invariant sets of scalar parabolic PDEs. We give several examples
of spatially inhomogeneous systems possessing infinite collections
of intricate stationary and time-periodic solutions.
\end{abstract}

\section{Introduction}
\label{s:intro}

We consider the invariant dynamics of one-dimensional second order
parabolic equations of the type
\begin{equation}\label{paracont}
  u_t = f(x,u,u_x,u_{xx}),
\end{equation}
where $u$ is a scalar function of the variables $t\in\re$ (time)
and $x \in S^1=\re/\ell\zet$ (periodic boundary conditions in
space), and $f$ is a $C^1$-function of its arguments.

\subsection{Assumptions}
\label{s:assumptions}

The case of periodic boundary conditions in $x$ provides richer
dynamics in general than Neumann or Dirichlet boundary conditions;
however, the techniques we introduce are applicable to a
surprisingly large variety of nonlinear boundary conditions.

This paper does not deal with the initial value problem, but
rather with the {\it bounded invariant dynamics}: bounded
solutions of \rmref{paracont} that exist for all time $t$. One
distinguishes three types of behaviors which are the building
blocks of all bounded invariant solutions to \rmref{paracont}
\cite{AngFied,Hale,Zel}.
\begin{enumerate}
\item[(i)]
{\em stationary patterns:}
$u(t,x) = u(x)$, $\forall t \in \re$,
\item[(ii)]
{\em periodic motions:}
$u(t+T,x) = u(t,x)$, for some period $T>0$,
\item[(iii)]
{\em homoclinic/heteroclinic connections:}
$\lim_{t\to\pm\infty} u(t,x) = u_\pm(x)$, where $u_\pm$
are stationary or periodic solutions of \rmref{paracont}.
\end{enumerate}

For the remainder of this paper we impose two natural assumptions
on \rmref{paracont}. The first is {\em uniform parabolicity}:
\vskip.2cm

\noindent
{\bf (f1)}\quad
$0<\lambda \le \partial_{w} f(x,u,v,w) \le \lambda^{-1}\,$,
uniformly $\,\forall\,(x,u,v,w) \in S^1\times \re^3$.
\vskip.2cm

\noindent This condition --- that \rmref{paracont} grows linearly
in $u_{xx}$ --- can be relaxed to degenerate parabolic equations
where the dependence on $u_{xx}$ is as a power law, see
\S\ref{sec_conc}. The second hypothesis is a {\em sub-quadratic
growth} condition on the $u_x$ term of $f$: \vskip.2cm

\noindent
{\bf (f2)}
\quad There exist constants $C>0$ such that

\quad $|f(x,u,v,w)| \le C(1+|v|^{\gamma})$, uniformly in both
$x\in S^1$ and on compact intervals in $u$ and $w$, for some
$0<\gamma<2$, \vskip.2cm

\noindent This will be necessary for regularity and control of
derivatives of solution curves, cf. \cite{AngFied}. This condition
is sharp: one can find examples of $f$ with quadratic growth in
$u_x$ for which solutions have singularities in $u_x$. Since our
topological data are drawn from graphs of $u$, the bounds on $u$
need to imply bounds on $u_x$ and $u_{xx}$: {\bf (f2)} does just
that.

A third {\em gradient} hypothesis will sometimes be assumed:

\noindent
{\bf (f3)}
\quad $f$ is {\em exact}, i.e.,
\begin{equation}\label{wexact}
    f(x,u,u_x,u_{xx}) = a(x,u,u_x)
    \Bigl[{d\over dx} \partial_{u_x} L - \partial_u L \Bigr],
\end{equation}
for a strictly positive and bounded function $a=a(x,u,u_x)$ and some
Lagrangian $L= L(x,u,u_x)$ satisfying
$0<\lambda \le a(x,u,u_x) \cdot \partial^2_{u_x} L(x,u,u_x) \le \lambda^{-1}$.
\vskip.2cm

In this case, we have a gradient system whose stationary solutions
are critical points of the action $\int L(x,u,u_x)dx$ over loops
of integer period in $x$. This condition holds for a wide variety
of systems. In general, systems with Neumann or Dirichlet boundary
conditions admit a {\em gradient-like} structure: there exists a
Lyapunov function which decreases strictly in $t$ along
non-stationary orbits. This precludes the existence of
nonstationary time-periodic solutions. It was shown by Zelenyak
\cite{Zel} that this gradient-like hypothesis holds for many
nonlinear boundary conditions which are a mixture of Dirichlet and
Neumann.

\subsection{Lifting the comparison principle}
\label{s:comp}

An important property of one-dimensional parabolic dynamics is the
lap-number principle of Sturm, Matano, and Angenent \cite{Ang1,
Matano1, Sturm} which, roughly, states that the number of nodal
regions in $x$ of $u(t,x)$ is a weak Lyapunov function for
\rmref{paracont}.

The lifting of this principle to the simultaneous evolution of
pairs of solutions is extremely fruitful. Consider two solutions
$u^1(t,x)$ and $u^2(t,x)$. Any tangency between the graphs
$u^1(t,\cdot)$ and $u^2(t,\cdot)$ at time $t=t_*$ is removed for
$t=t_*+\epsilon$ (for all small $\epsilon>0$) so as to {\em
strictly} decrease the number of intersections of the graphs. This
holds even for highly degenerate tangencies of curves \cite{Ang1}.
As shown in the work of Fiedler and Mallet-Paret \cite{FMP}, this
comparison principle implies that the dynamics of \rmref{paracont}
is weakly Morse-Smale (all bounded orbits are either fixed points,
periodic orbits, or connecting orbits between these), see
\cite{Hale, Zel}.

The idea behind this paper, following the discrete version of this
phenomenon in \cite{GVV}, is to ``lift'' the comparison principle
from pairs of solutions to larger ensembles of solution curves.
The local data attached to pairs of curves --- intersection number
--- can be lifted to more global data about patterns of
intersections via the language of topological {\em braid theory}.
A similar theory for geodesics on two dimensional surfaces has
been developed in \cite{Ang2}, and has served as a guideline for
some of the ideas used here.

Consider a collection
$\u=\left\{u^\alpha(t,\cdot)\right\}_{\alpha=1}^n$ of $n>1$
solutions to \rmref{paracont}, where, to obey the periodic
boundary conditions in $x$,
$\left\{u^\alpha(t,0)\right\}_{\alpha=1}^n =
\left\{u^\alpha(t,1)\right\}_{\alpha=1}^n$ as sets of
points.\footnote{This condition permits solutions with integral
period which ``wrap'' around the circle.} Instead of thinking of
the graphs of $u^\alpha(t,\cdot)$ as being evolving curves in the
$(x,u)$ plane, we take the {\em 1-jet extension} of each curve and
think of it as an evolving curve in $(x,u,u_x)$ space.
Specifically, for each $t$, $u^\alpha(t,\cdot):[0,1]\to
[0,1]\times\real^2$ given by
$x\mapsto(x,u^\alpha(x,t),u^\alpha_x(x,t))$. As long as these
curves do not intersect in their 3-d representations, we have what
topologists call a {\em braid}. In particular, such a braid is
said to be {\em closed} (the ends $x=0$ and $x=1$ are identified)
and {\em Legendrian} (the curves are all tangent to the standard
contact structure $dx_2-x_3\,dx_1=0$).

As these curves evolve under the PDE, the topological type of the
braid can change. The topological equivalence class of a closed
Legendrian braid is the appropriate analogue of the intersection
data for pairs of curves. Indeed, there is a natural group
structure on braids with $n$ strands. We argue in a ``braid
theoretic'' version of the comparison principle that the {\em
algebraic length} of a braid given by solutions
$\{u^\alpha(t)\}_{\alpha=1}^n$ is a weak Lyapunov function for the
dynamics of \rmref{paracont}.

\subsection{Main results}
\label{s:results}

The goal of this paper, following earlier work in \cite{GVV} on a
discrete version of this problem, is to define an index for closed
Legendrian braids and to use this as the basis for detecting
invariant dynamics of \rmref{paracont}. See \S\ref{sec_braids} for
definitions and background on the discrete version.

For purposes of detecting invariant dynamics of \rmref{paracont},
we work with braids $\u$ relative to some fixed braid $\v$. One
thinks of $\v$ as a braid for which dynamical information is
known, namely, that its strands are $t$-invariant solutions to
\rmref{paracont}, the entire set of which respects the periodic
boundary conditions (individual strands might not: see
Fig.~\ref{fig_braids}). One thinks of $\u$ as consisting of
``free'' strands about which nothing is known with regards to
dynamical behavior.

We show that there exists a well-defined {\em homotopy index} that
maps a (closed, Legendrian, relative) braid class represented by
$\{\u\rel\v\}$ to a pointed homotopy class of spaces,
$\HH(\u\rel\v)$. This index is at heart a Conley index for a
suitable configuration space which is isolated thanks to the
braid-theoretic comparison principle. A coarser homology index
sends such braids to a polynomial $P_\tau(\HH)$ in one variable,
$\tau$.

The main results of this paper are forcing theorems for stationary
and periodic solutions.

\subsubsection{Stationary solutions}
For our main results we restrict to braid classes which have two
compactness properties: proper and bounded. Roughly speaking, a
relative braid class $\{\u\rel\v\}$ is {\em proper} if none of the
components of $\u$ can be collapsed onto $\u$ or $\v$. A relative
braid class $\{\u\rel\v\}$ is {\em bounded} if all strands of $\u$
are uniformly bounded with respect to all representatives
$\u\rel\v$ of the braid class: see \S\ref{sec_braids}.

\begin{theorem}
\label{thm_Main} Let \rmref{paracont} satisfy {\bf (f1)} and {\bf
(f2)} with $\v$ a stationary braid. If $\{\u\rel\v\}$ is a bounded
proper braid class, then there exists a stationary solution of
this braid class if the Euler characteristic of $\HH$,
$\chi(\HH):=P_{-1}(\HH)$, is nonvanishing. If in addition $f$
satisfies {\bf (f3)}, then there are at least $\abs{P_\tau(\HH)}$
stationary solutions of this braid class, where $\abs{\cdot}$
denotes the number of nonzero monomials.
\end{theorem}

The above theorem is formulated for periodic boundary conditions.
In the case of other boundary conditions the Zelenyak result
implies that \rmref{paracont} is automatically gradient-like so
that the second part of Theorem \ref{thm_Main} is superfluous is
those cases.

\begin{remark}{\em
\label{rem_nondeg} With additional knowledge, $P_\tau(\HH)$ can
reveal more of the dynamics. For example, assume for simplicity
that the invariant sets are known to be hyperbolic and that the
strands of $\u$ form a single-component braid (the graph of $\u$
is connected as a subset of $S^1\times\real$). In this setting,
the strong Morse inequalities yield more information on
multiplicity of solutions. As pointed out before the critical
elements of \rmref{paracont} are equilibrium solutions and
periodic orbits. Therefore the Morse relations are given by
$$
\sum_i a_i \tau^i + \sum_j b_j\tau^j(1+\tau) = P_\tau(\HH) +
(1+\tau) Q_\tau,
$$
where $Q_\tau $ is a polynomial with non-negative coefficients.
The coefficients $a_i$ count the number of equilibrium solutions
of Morse index $i$, while the $b_j$ count the number of periodic
orbits of Morse index $j$. If one assumes nondegeneracy, then the
Morse relations can be used to compute $P_\tau(\HH)$.
}\end{remark}

\begin{remark}{\em
\label{rem_infinite} In the exact case the lower bound on the
number of critical points can refined even further. For parabolic
recurrence relations the spectrum of a critical point satisfies
$\lambda_0<\lambda_1\le \lambda_2 <\lambda_3\le
\lambda_4<\lambda_5 ....$. This ordering has special bearing on
non-degenerate critical points with odd index. To be more precise,
for a `topological' non-degenerate critical point $\u$ with
$P_\tau(\u) = A\tau^{2k+1}$ it holds that $A=1$. More details of
this can be found in \S\ref{sec_Main}. A direct consequence there
are at least as many critical points as the sum of the odd Betti
numbers of $\HH$. If we write $P_\tau(\HH) = P^{\rm
odd}_\tau(\HH)+P^{\rm even}_\tau(\HH)$, then our lower bound on
the number of critical points becomes
$$
P^{\rm odd}_1(\HH) + |P^{\rm even}_\tau(\HH)|,
$$
which lies in between $|P_\tau(\HH)|$ and $P_1(\HH)$.
}\end{remark}

The proof of Theorem~\ref{thm_Main} appears in \S\ref{sec_proof}.
First, however, we introduce the relevant portions of braid theory
(\S\ref{sec_braids}), followed by a review
(\S\ref{sec_inv}-\ref{sec_dyn}) of the discrete braid index
constructed in \cite{GVV}.

This theory applies to a wide array of inhomogeneous equations. In
\S\ref{sec_exstat} we show,

\begin{example}{\em
The equation
    \begin{equation}
    \label{eq_Stat}
       u_t = u_{xx} - \frac{5}{8}\sin 2x\, u_x
       + \frac{\cos x}{\cos x + \frac{3}{\sqrt{5}}}
       u(u^2-1).
    \end{equation}
possesses stationary solutions in an infinite number of distinct
braid classes. As a matter of fact we show that one can embed an
Bernoulli shift into the stationary equation. }
\end{example}

\begin{example}{\em
For $\epsilon\ll 1$ and any smooth nonconstant $h:S^1\to(0,1)$,
the equation
    \begin{equation}
    \label{eq_Nakashima}
        \epsilon^2 u_t = \epsilon^2 u_{xx} + h(x) u (1-u^2) .
    \end{equation}
possesses stationary solutions spanning an infinite collection of
braid classes. This example was studied by Nakashima
\cite{Nakashima1,Nakashima2}.
}\end{example}

These two examples can be generalized greatly.
Theorem~\ref{thm_infinite} gives extremely broad conditions which
force an infinite collection of stationary solutions.

\subsubsection{Periodic solutions}

We also lay the foundation for using the braid index to find
time-periodic solutions. For simplicity in the analysis, we
restrict our attention to equations of the form
\begin{equation}
\label{paracontII}
    u_t = u_{xx} + g(x,u,u_x) ,
\end{equation}
which trivially satisfies Hypothesis {\bf (f1)}. By assuming
Hypothesis {\bf (f2)} (without the $w$ variable), we prove an
analogue of Theorem~\ref{thm_Main} for time-periodic solutions of
\rmref{paracontII}. As we pointed out before, time-periodic
solutions can exist by the grace of the boundary conditions. As
the result of Zelenyak implies, in most cases a weak version of
{\bf (f3)} holds (gradient-like) and the only critical elements
are stationary solutions.

\begin{remark}{\em
\label{rem_rotating} A fundamental class of time-periodic orbits
are the so-called {\em rotating waves}. For an equation which is
autonomous in $x$, one makes the rotating wave hypothesis that
$u(t,x) = U(x-ct)$, where $c$ is the unknown wave speed.
Stationary solutions for the resulting equation on $U(\xi)$ yield
rotating waves. Modulo the unknown wave speed --- a nonlinear
eigenvalue problem --- Theorem \ref{thm_Main} now applies. In
\cite{AngFied} it was proved that time-periodic solutions are
necessarily rotating waves for an equation autonomous in $x$.
However, in the non-autonomous case, the rotating wave assumption
is highly restrictive.
}\end{remark}

We present a very general technique for finding time-periodic
solutions without the rotating wave hypothesis.

\begin{theorem}
\label{thm_Main2} Let \rmref{paracontII} satisfy {\bf (f2)} with
$\v$ a stationary braid. Let $\{\u\rel\v\}$ be a bounded proper
braid class with $\u$ a single-component braid and
$P_\tau(\HH)\neq 0$. If the braid class is not stationary for
\rmref{paracontII}, then there exists a time-periodic solution in
this braid class.
\end{theorem}

\begin{remark}{\em
\label{rem_nec} In certain examples one can find braid classes in
which a given equation cannot have stationary solutions. Since the
only possible critical elements in that case are periodic orbits
it follows that the Poincar\'e polynomial   has to be of the form
$P_\tau(\HH) = (1+\tau) p_\tau(\HH)$. The polynomial $p_\tau(\HH)$
gives a lower bound on the number of periodic orbits (in the
non-degenerate case). The single-component hypothesis on $\u$
(namely, that the graph of $\u$ is connected in $S^1\times\real$)
is not crucial. For free strands forming multi-component braids
$\u$, each component of $\u$ will be time-periodic. Their periods
may not be rationally related, however, leading to a
quasi-periodic solution in time in the multi-component braid
class. }
\end{remark}

It was shown in \cite{AngFied} that a singularly perturbed van der
Pol equation,
\[
    u_t = \epsilon u_{xx} + u(1-\delta^2u^2) + u_xu^2,
\]
possesses an arbitrarily large number of rotating waves depending
on $\epsilon\ll 1$ for fixed $0<\delta$. We generalize their
result:

\begin{example}{\em
Consider the equation
    \begin{equation}
        u_t = u_{xx}+ug(u)+u_xh(x,u,u_x) ,
    \end{equation}
where the non-linearity is assumed to satisfy {\bf (f2)}, i.e. $h$
has sub-linear growth in $u_x$ at infinity. Moreover, $g$ and $h$
satisfy the following hypotheses:
\begin{itemize}
\item[{\bf(g1)}]
    $g(0)>0$, and $g$ has at least one positive and one negative root;
\item[\bf{(g2)}]
    $h>0$ on $\{ u u_x\neq 0 \}$.
\end{itemize}
Then this equation possesses time-periodic solutions spanning an
infinite collection of braid classes.
}\end{example}

We provide details in \S\ref{sec_exper}. All of the periodic
solutions implied are dynamically unstable. In the most general
case (those systems with $x$-dependence), the periodic solutions
are not rigid rotating waves.

\section{Braids}
\label{sec_braids}

The results of this paper require very little of the extensive
theory of braids developed by topologists \cite{Bir}. However,
since the definitions motivate our constructions, we give a brief
tour.

\subsection{Topological braids}

A topological braid on $n$ strands is an embedding
$\beta:\coprod_1^n[0,1]\hookrightarrow\real^3$ of a disjoint union
of $n$ copies of $[0,1]$ into $\real^3$ such that
\begin{itemize}
\item[(a)] the left endpoints $\beta(\coprod_1^n\{0\})$ are
$\{(0,i,0)\}_{i=1}^n$;
\item[(b)] the right endpoints $\beta(\coprod_1^n\{1\})$ are
$\{(1,i,0)\}_{i=1}^n$; and
\item[(c)] $\beta$ is transverse to the planes $x_1={\rm constant}$.
\end{itemize}
Two braids are said to be of the same {\em topological braid
class} if they are homotopic in the space of braids: one braid
deforms to the other without any intersections of the strands. A
{\it closed} topological braid is obtained if one quotients out
the range of the braid embeddings via the equivalence relation
$(0,x_2,x_3)\sim(1,x_2,x_3)$ and alters the restrictions (a) and
(b) of the position of the endpoints to be
$\beta(\coprod_i^n\{0\})=\beta(\coprod_1^n\{1\})$. Thus, a closed
braid is a collection of disjoint embedded arcs in
$[0,1]\times\re^2$ (with periodic boundary conditions in the first
variable) which are everywhere transverse to the planes
$x_1=$constant.

In this paper, we restrict attention to those braids whose strands
are of the form $(x,u(x),u_x(x))$ for $0\leq x\leq 1$. These are
sometimes called {\em Legendrian braids} as they are tangent to
the canonical contact structure $dx_2-x_3\,dx_1$. No knowledge of
Legendrian braid theory is assumed for the remainder of this work,
but we will use the term freely to denote those braids lifted from
graphs.

\subsection{Braid diagrams}
\label{s:diagrams}

The specification of a topological braid class (closed or
otherwise) may be accomplished unambiguously by a labeled
projection to the $(x_1,x_2)$-plane; a {\it braid diagram}.
Labeling is done as follows: perturb the projected curves slightly
so that all strand crossings in the projection are transversal and
disjoint. Then, mark each crossing via $(+)$ or $(-)$ to indicate
whether the crossing is ``left over right'' or ``right over left''
respectively.

Since a Legendrian braid is of the form $(x,u(x),u_x(x))$, no such
marking of crossings in the $(x,u)$ projection are necessary: all
crossings have positive labels. For the remainder of this paper we
will consider only such {\em positive braid diagrams}. We will
analyze parabolic PDEs by working on spaces of such braid
diagrams. Although Legendrian braids are the right types of braids
to work with as solutions to \rmref{paracont} (cf. the smoothing
of initial data for heat flow), our discretization techniques will
require a more robust $C^0$ theory for braid diagrams. Thus, we
work on spaces of braid diagrams with topologically transverse
strands:

\begin{definition}\label{PB}
The space of {\em closed positive braid diagrams on $n$ strands},
denoted $\Conf^n$, is the space of all pairs $(\uu,\tau)$ where
$\tau\in S_n$ is a permutation on $n$ elements, and
$\uu=\{u^\alpha(x)\}_{\alpha=1}^n$ is an unordered collection of
$H^1$-functions --- {\em strands} --- satisfying the following
conditions:
\begin{enumerate}
\item[(a)] {\bf Periodicity:}
$u^\alpha(1) = u^{\tau(\alpha)}(0)$ for all $\alpha$.
\item[(b)] {\bf Transversality:}
for any $\alpha\neq\alpha'$ such that $u^\alpha(x_*)=u^{\alpha'}(x_*)$
for some $x_*\in [0,1]$, it holds that $u^\alpha(x)-u^{\alpha'}(x)$
has an isolated sign change at $x=x_*$.
\end{enumerate}
Because the strands of $\uu$ are unordered, we naturally identify
all pairs $(\uu,\tau)$ and $(\uu,\tilde\tau)$ satisfying
$\tilde\tau=\sigma\tau\sigma^{-1}$ for some permutation $\sigma\in
S_n$. Henceforth we suppress the permutations $\tau$ from the
description of a braid, it being understood implicitly.
\end{definition}

The path components of $\Conf^n$ comprise the {\it braid classes}
of closed positive braid diagrams. The braid class of a braid
diagram $\u$ is denoted by $\{\u\}$. Any braid diagram $\u$ with
$C^1$-strands naturally lifts to a Legendrian braid by the 1-jet
extension of $u^\alpha$ to the curve
$(x,u^\alpha(x),u_x^\alpha(x))$. If we allow the strands to
intersect --- disregarding condition (b) of Definition~\ref{PB}
--- we obtain a closure of the space $\Conf^n$, which we denote
$\overline{\Conf^n}$. The `discriminant' $\Sigma^n :=
\overline{\Conf^n} - \Conf^n$ defines the {\it singular} braid
diagrams.

\subsection{Discrete braid diagrams}

From topological braids we have passed to braid diagrams in order
to describe invariant curves for parabolic PDEs. There is one last
transformation we must impose: a spatial discretization.

\begin{definition}\label{PL}
The space of {\em period $d$ discrete braid diagram on $n$
strands}, denoted $\DConf^n_d$, is the space of all pairs
$(\uu,\tau)$ where $\tau\in S_n$ is a permutation on $n$ elements,
and $\uu=\{u^\alpha\}_{\alpha=1}^n$ is an unordered collection of
vectors $u^\alpha=(u^\alpha_i)_{i=0}^d$ --- {\em strands} ---
satisfying the following conditions:
\begin{enumerate}
\item[(a)] {\bf Periodicity:}
    $u^\alpha_d = u^{\tau(\alpha)}_0$ for all $\alpha$.
\item[(b)] {\bf Transversality:}
    for any $\alpha\neq\alpha'$
    such that $u^\alpha_i=u^{\alpha'}_i$ for some $i$,
\begin{equation}
\label{eq_transverse}
        \bigl(u^\alpha_{i-1}-u^{\alpha'}_{i-1}\bigr)
    \bigl(u^\alpha_{i+1}-u^{\alpha'}_{i+1}\bigr) < 0.
\end{equation}
\end{enumerate}
As in Definition~\ref{PB}, the permutation $\tau$ is defined up to
conjugacy (since the strands are unordered) and will henceforth
not be explicitly written.
\end{definition}

The path components of $\DConf^n_d$ comprise the {\it discrete
braid classes} of period $d$. The discrete braid class of a
discrete braid diagram $\u$ is denoted $[\u]$.
If we disregarding condition (b) of Definition~\ref{PL}, we obtain a closure
of the space $\DConf^n_d$, which we denote $\overline{\DConf^n_d}$.
The `discriminant' $\Sigma^n_d:= \overline{\DConf^n_d} - \DConf^n_d$ defines
the {\it singular} discrete braid diagrams of period $d$.

\begin{figure}[hbt]
\begin{center}
\psfragscanon
\psfrag{x}[][]{$x$}
\psfrag{u}[r][]{$u$}
\psfrag{d}[tr][]{$u_x$}
\includegraphics[width=5.25in]{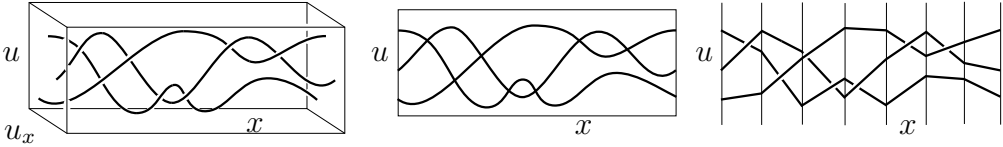}
\caption{Three types of braids: a Legendrian topological braid [left],
its braid diagram [center], and a discrete braid diagram
[right].}
\label{fig_braids}
\end{center}
\end{figure}

Figure~\ref{fig_braids} summarizes the three types of braids
introduced in this section.

\subsection{Discretization: back and forth}

It is straightforward to pass from topological to discrete braids
and back again.

\begin{definition}
\label{def_disc} Let $\u\in\Conf^n$ be a topological closed braid
diagram. The {\em period-$d$ discretization} of $\uu$ is defined
to be
\begin{equation}
\disc_d(\uu) = \{\disc_d(u^\alpha)\}^\alpha := \{u^\alpha(i/d)\}_i^\alpha.
\end{equation}
Conversely, given a discrete braid $\uu\in\DConf^n_d$, we
construct a piecewise-linear [PL] topological braid diagram,
$\pl(\uu):=\{\pl(u^\alpha)\}$, where $\pl(u^\alpha)$ is the
$C^0$-strand given by
\begin{equation}
\label{eq_PL}
   \pl(u^\alpha)(x) := u_{\lfloor d\cdot x\rfloor} +
   \bigl(d\cdot x - \lfloor d\cdot x\rfloor\bigr)
   \bigl( u_{\lceil d\cdot x\rceil} - u_{\lfloor d\cdot x\rfloor}\bigr) .
\end{equation}
\end{definition}

The following lemma is left as an exercise.
\begin{lemma}
Let $\uu\in\Conf^n$ and $\vv\in\DConf^n_d$.
\begin{enumerate}
\item
  $\pl$ sends the discrete braid class $[\vv]$ to a well-defined
  topological braid class $\{\pl(\vv)\}$.
\item
  For $d$ sufficiently large, $\{\pl(\disc_d(\uu))\}=\{\uu\}$.
\end{enumerate}
\end{lemma}

The second part of this lemma accommodates the obvious fact that
braiding data is lost if the discretization is too coarse. This
leads to the following definition:
\begin{definition}\label{disc}
A discretization period $d$ is {\em admissible} for $\u \in
\Conf^n$ if
\[ \{\pl(\disc_d(\u))\}=\{\u\} . \]
\end{definition}

In the next section, we will describe a Morse-Conley topological
index for pairs of braids which relies on algebraic length of the
braid as a Morse function. Rather than detail the algebraic
structures, we use an equivalent geometric formulation of length:

\begin{definition}
\label{def_length} The length of a topological braid
$\uu\in\Conf^n$, denoted $\length(\uu)$, is defined to be the
total number of intersections in the braid diagram. If
$\uu\in\DConf^n_d$ is a discrete braid, then
$\length(\uu):=\length(\pl(\uu))$.
\end{definition}

\section{Braid invariants}
\label{sec_inv}

We give a concise description of the invariant of \cite{GVV} for
relative discrete closed braids.

\subsection{Relative braids}

The motivation for the homotopy braid index is a forcing theory:
given a stationary braid $\vv$, does it force some other braid
$\uu$ to also be stationary with respect to the dynamics? This
necessitates understanding how the strands of $\uu$ braid relative
to those of $\vv$.

\begin{definition}
Given $\v \in \Conf^m$, define
\[
    \Conf^n{\rel}\v := \{ \u \in \Conf^n~:~\u\cup\v \in \Conf^{n+m}\} .
\]
The path components of $\Conf^{n}\rel\v$, comprise the {\em
relative braid classes}, denoted $\{\u\rel\v\}$. In this setting,
the braid $\v$ is called the {\em skeleton}.
\end{definition}

This procedure partitions $\Conf^n$ relative to $\v$: not only are
tangencies between strands of $\uu$ illegal, so are tangencies
with the strands of $\vv$.

The definitions for discrete relative braids are analogous.
\begin{definition}
Given $\v \in \DConf_d^m$, define
\[
    \DConf_d^n{\rel}\v := \{\u\in\DConf_d^n~:~\u\cup\v \in \DConf_d^{n+m}\} .
\]
The path components of $\DConf_d^{n}\rel\v$, comprise the {\em
relative discrete braid classes}, denoted $[\u\rel\v]$. In this
setting, the braid $\v$ is called the {\em skeleton}.
\end{definition}

The operations $\disc_d$ and $\pl$ have obvious extensions to
relative braids by acting on both $\uu$ and $\vv$.

\subsection{Bounded and proper relative braids}

\begin{definition}\label{proper}
A relative braid class $\{\u\rel\v\}$ is called {\em proper} if it
is impossible to find an isotopy $\u(t)\rel\v$ such that $\u(0) =
\u$, $\u(t)\rel\v \in \{\u\rel\v\}$, for $t\in [0,1)$, and
$\u(1)\cup\v \in \Sigma^{n+m}$ is a diagram where an entire
component of the braid $\u(1)$ has collapsed onto itself, another
component of $\u(1)$, or a component of $\v$. A discrete relative
braid class is proper if it is the discretization of a proper
topological relative braid class.
\end{definition}

The index we define is based on the topology of a relative braid
class. It is most convenient to define this on compact spaces;
hence the following definition.

\begin{definition}
A braid class (topological or discrete) is {\em bounded} if
$\{\u\rel\v\}$ is a bounded set in $\overline{\Conf^n}$.
\end{definition}

\begin{figure}[hbt]
\begin{center}
\includegraphics[width=4in]{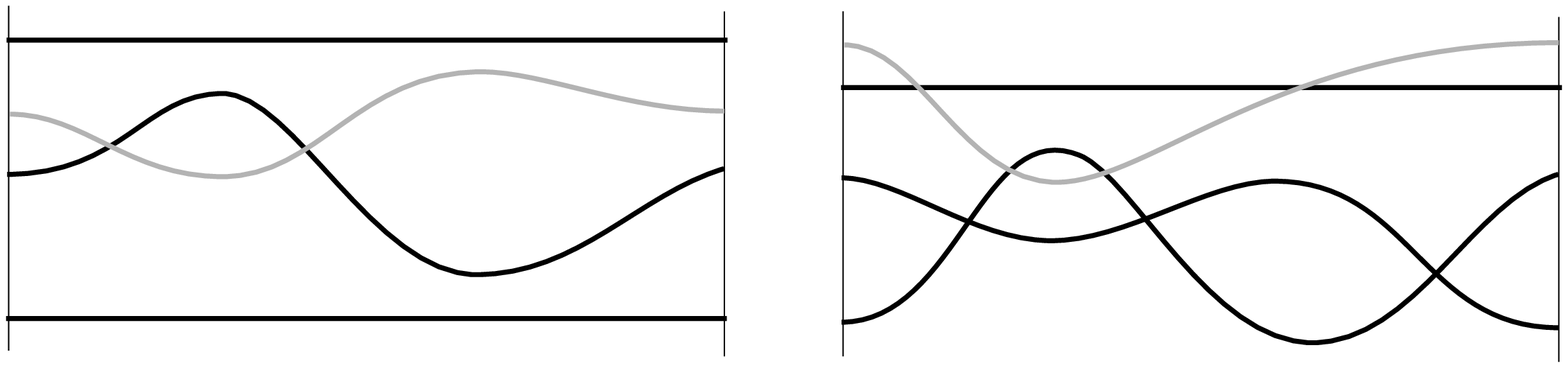}
\caption{[left] a bounded but improper braid class ; [right] a proper,
but unbounded braid class. Black strands are fixed, grey free.}
\label{numfig3}
\end{center}
\end{figure}

For the remainder of the paper, all braids will be assumed proper
and bounded unless otherwise stated.

\subsection{The Conley index for braids}
\label{sec_CI}

Consider a discrete relative braid class
$[\uu\rel\vv]\subset\DConf^n_d$ which is bounded and proper. We
associate to this class a Conley-type index for a class of
dynamics on spaces of discrete braids. This will become an
invariant of topological braids via discretization.

Denote by $N$ the closure of $[\uu\rel\vv]$ in the space
$\overline{\DConf^n_d}\rel\vv$. We identify an ``exit set'' on the
boundary of $N$ consisting of those relative braids whose length
$\length$ can be decreased by a small perturbation. Let
$\ww\in\partial N$ denote a singular braid on the boundary of $N$
and let $W$ be a sufficiently small neighborhood of $\ww$ in
$\overline{\DConf^n_d}\rel\vv$. Then $W$ is sliced by
$\Sigma^n_d\rel\vv$ into a finite number of connected components
representing distinct neighboring braid classes, each component
having a well-defined braid length $\iota\in\zed^+$. Define the
{\em exit set}, $N^-$, of $N$ to be those singular braids at which
$\length$ can decrease:
\begin{equation}
\label{eq_exit}
  N^- := \cl\left\{ \ww\in\partial N : \length
                  {\mbox{ is locally maximal on {\sc int}}}(N) \right\} ,
\end{equation}
where $\cl$ denotes closure in $\partial N$
\begin{definition}
The {\em Conley index} of a discrete (proper, bounded, relative)
braid class $[\uu\rel\vv]$ is defined to be the pointed homotopy
class of spaces
\begin{equation}
\label{eq_h}
  h([\uu\rel\vv]) = \left[ N/N^- \right] := \left(N/N^-, [N^-]\right) .
\end{equation}
\end{definition}

\begin{example}{\em
\label{ex_indexex} Consider the period-2 braid illustrated in
Fig.~\ref{fig_indexex}[left] possessing exactly one free strand
with anchor points $u_1$ and $u_2$. The anchor point in the
middle, $u_1$, is free to move vertically between the fixed points
on the skeleton. At the endpoints, one has a singular braid in
$\Sigma$ which is on the exit set since a slight perturbation
sends this singular braid to a different braid class with fewer
crossings. The end anchor point, $u_2$, can move vertically
between the two fixed points on the skeleton. The singular
boundaries are in this case {\it not} on the exit set since
pushing $u_2$ across the skeleton increases the number of
crossings.

\begin{figure}[hbt]
\begin{center}
\psfragscanon
\psfrag{0}[][]{}
\psfrag{1}[][]{}
\psfrag{2}[][]{}
\psfrag{a}[lt][]{$u_2$}
\psfrag{b}[][]{$u_1$}
\includegraphics[width=5.25in]{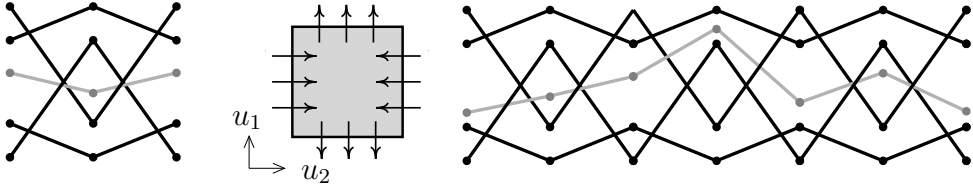}
\caption{A period two braid [left], the associated configuration
space [center], and a period six generalization [right].}
\label{fig_indexex}
\end{center}
\end{figure}

Since the points $u_1$ and $u_2$ can be moved independently, the
configuration space $N$ in this case is the product of two compact
intervals. The exit set $N^-$ consists of those points on
$\partial N$ for which $u_1$ is a boundary point. Thus, the
homotopy index of this relative braid is $[N/N^-]\simeq S^1$.

By taking a chain of copies of this skeleton (i.e., taking a cover
of the spatial domain), one can construct examples with one free
strand weaving in and out of the fixed strands in such a way as to
produce an index with homotopy type $S^k$ for any $k\geq 0$.
}\end{example}

The extension of the Conley index to topological braid diagrams is
straightforward: choose an admissible discretization period $d$,
take the Conley index of the period-$d$ discretization, then show
that this is independent of $d$. The key step --- independence
with respect to $d$ --- is, unfortunately not true. For $d$
sufficiently small, there may be different discrete braid classes
which define the same topological braid. The information from any
one of these coarse components is incomplete. The following
theorem, which is the main result from \cite{GVV}, resolves this
obstruction.

\begin{theorem}[see \cite{GVV}, Thm. 19 and Prop. 27]
\label{thm_GVV1} For $d$ sufficiently large,\footnote{A sufficient
though high lower bound is the number of crossings of $\u$ with
itself and with $\v$.} the Conley index
$h([\disc_d\u\rel\disc_d\v])$ is independent of $d$ and thus an
invariant of the topological braid class $\{\u\rel\v\}$.
\end{theorem}

\begin{definition}
\label{defH} Given a topological braid class $\{\u\rel\v\}$,
define the {\em homotopy index} to be
\begin{equation}
\label{eq_HH}
     \HH(\u\rel\v) := h([\disc_d\u\rel\disc_d\v]) .
\end{equation}
for $d$ sufficiently large.
\end{definition}

For purposes of this paper, the homotopy index is defined with $d$
sufficiently large. This is well-defined, but not optimal for
doing computations. To that end, one can use the more refined
formula of \cite{GVV}, which computes $\HH$ for any admissible
discretization period $d$ via wedge sums: we will not require this
complication in this paper.

For most applications it suffices to use the homological
information of the index given by its Poincar\'e polynomial
\begin{equation}
  P_\tau(\HH) := \sum_{k=0}^\infty {\rm dim~} H_k(\HH)\tau^k
  = \sum_{k=0}^\infty {\rm dim~} H_k(N,N^-)\tau^k .
\end{equation}
This also has the pleasant corollary of making the index
computable via rigorous homology algorithms.

\section{Dynamics and the braid index}
\label{sec_dyn}

The homotopy braid invariant is defined as a ``Conley index.''
This index has significant dynamical content.

The most basic version of the Conley index has the following
ingredients \cite{Conley}: given a continuous flow on a metric
space, a subset $N$ is said to be an {\em isolating block} if all
points on $\partial N$ leave $N$ under the flow in forwards and/or
backwards time. The Conley index of $N$ with respect to the flow
is then the pointed homotopy class $[N/N^-]$, where $N^-$ denotes
the {\em exit set}, or points on $\partial N$ which leave $N$
under the flow in forwards time. Standard facts about the index
include (1) invariance of the index under continuous changes of
the flow and the isolating block; and (2) the forcing result: a
nonzero index implies that the flow has an invariant set in the
interior of $N$. In order to implement Conley index theory in
combination with braids we define the following class of dynamical
systems.

\begin{definition}
\label{def_PRR} Given $d>0$, a {\em parabolic recurrence relation}
$\R$ on $\zed/d\zed$ is a collection of $C^1$-functions
$\R_i:\real^3\to\real$, $i\in\zed/d\zed$ such that for each $i$,
$\del_1\R_i>0$ and $\del_3\R_i\geq 0$. We say that $\R$ is {\em
exact} if there exists a sequence of $C^2$-generating functions
$S_i$ such that
\begin{equation}
    \R_i(u_{i-1},u_i,u_{i+1}) =
    \partial_2S_{i-1}(u_{i-1},u_i) +
    \partial_1S_i(u_i,u_{i+1}) \quad \forall i .
\end{equation}
\end{definition}

A parabolic recurrence relation (henceforth {\em PRR}) defines a
vector field on $\overline{\DConf^n_d}$,
\begin{equation}
\label{eq_PRR}
  \frac{d}{dt}(u_i^\alpha) = \R_i(u_{i-1}^\alpha,u_i^\alpha,u_{i+1}^\alpha) ,
\end{equation}
with all subscript operations interpreted modulo the permutation
$\tau$: $u_{d+1}^\alpha=u_1^{\tau(\alpha)}$. The flow generated by
\rmref{eq_PRR}  is called a {\em parabolic flow} on $\DConf_d^n$.
For more details see \cite{GVV}. Exact PRR's induce a flow which
is the gradient flow of
$W(\u):=\sum_iS_i(u^\alpha_i,u^\alpha_{i+1})$.

\begin{figure}[hbt]
\begin{center}
\psfragscanon
\psfrag{u}[rb][]{\large $\tilde\uu$}
\psfrag{S}[t][]{\large $\Sigma$}
\psfrag{0}[][]{\large $i-1$}
\psfrag{1}[][]{\large $i$}
\psfrag{B}[][]{\large $[\uu\rel\vv]$}
\psfrag{2}[][]{\large $i+1$}
\includegraphics[angle=0,width=5.25in]{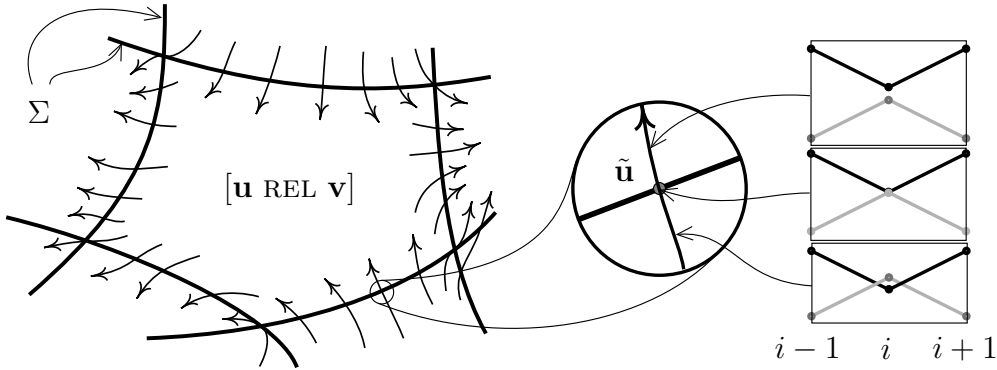}
\caption{A parabolic flow on a (bounded and proper) braid class is
transverse to the boundary faces, making the braid class into an
isolating block. The local linking of strands decreases strictly
along the flow lines at a singular braid $\tilde\uu$. }
\label{fig_transverse}
\end{center}
\end{figure}

A  parabolic flow acts on discrete braid diagrams in much the same
way that \rmref{paracont} acts on topological braid diagrams. As
we have defined it in \S\ref{sec_CI}, the Conley index for a
discrete braid class $[\u\rel\v]$ uses its closure $N={\rm
cl}[\u\rel\v]$ as an isolating block. Indeed, if $[\u\rel\v]$ is a
bounded, then $N$ is a compact set. If $[\u\rel\v]$ is proper,
then vector field on $\overline{\DConf^n_d }\rel \v$ induced by
$\R$, is transverse to $\partial N$, and $N$ is really an
isolating block for the parabolic flow. The set $N^-$ defined in
the previous section then  is the exit for $N$. This particular
link lies at the heart of the theory and follows from the a
discrete version of the comparison principle \cite{FOl, MPSm,
Smillie}. Details of the construction can be found in \cite{GVV},
where it is shown that the index $h([\u\rel\v])$ defined via Eqns.
(\ref{eq_exit}) and (\ref{eq_h}) is the Conley index of any PRR
which fixes $\v$. Fig.~\ref{fig_transverse} illustrates the action
of a parabolic flow on braids.

In \cite{GVV} it furthermore is shown that certain Morse
inequalities hold for stationary solutions of \rmref{eq_PRR}. The
Morse inequalities also provide information about the periodic
orbits. This is due to the fact that for parabolic systems the set
of bounded solutions consists only of stationary points, periodic
orbits, and connections between them.

Theorem~\ref{thm_Main} is an extension of the following results
for parabolic lattice systems.

\begin{theorem}\cite{GVV}
\label{thm_GVV2} Let $\R$ be a parabolic recurrence relation. The
induced flow on a bounded proper discrete braid class
$[\u\rel\v]$, where $\v$ is a stationary skeleton, has an
invariant solution within the class $[\u\rel\v]$ if the Conley
index $h = h([\u\rel\v])$ is nonzero. Furthermore:
\begin{enumerate}
\item
    If the Euler characteristic $\chi(h) \not = 0$ then
        there exist   stationary solutions of braid class
    $[\u\rel\v]$.
\item
    If $\R$ is exact, then the number of stationary solutions
    of braid class $[\u\rel\v]$ is bounded below by
    $\abs{P_\tau(h)}$, the number of nonzero monomials of the
    Poincar\'e polynomial of the index.
\end{enumerate}
\end{theorem}

If a proper bounded braid class $[\u\rel \v]$ contains no
stationary braids for a particular recurrence relation $\R$, then
$h(\u\rel\v) \not =0$ forces periodic solutions of \rmref{eq_PRR},
i.e. the components of $\u$ are periodic. If the system is
non-degenerate the number of orbits is given by $P_1(h)/2$. As a
consequence in this case $P_\tau(h)$ is divisible by $1+\tau$ and
$\R$ is {\em not} exact. Note that for $d$ large enough the
topological information is contained in the invariant $\HH$ for
the topological braid class $\{\u\rel\v\}$.

\section{Examples: stationary solutions}
\label{sec_exstat}

The following examples all satisfy Hypotheses {\bf (f1)} and {\bf
(f2)}.

\begin{example}
\label{ex_Matano}{\em
Consider the following family of spatially inhomogeneous
Allen-Cahn equations studied by Nakashima
\cite{Nakashima1,Nakashima2}:
\begin{equation}
\epsilon^2 u_t = \epsilon^2 u_{xx} + h(x) u (1-u^2) ,
\end{equation}
where $h:S^1\to(0,1)$ is not a constant. Clearly this equation has
stationary solutions $u=0,\pm 1$ and is exact with Lagrangian
\[
    L = \frac{1}{2}\epsilon^2 u_x^2
    - \frac{1}{4}\left(h(x)u^2(2-u^2)\right) .
\]
According to \cite{Nakashima1}, for any $N>0$, there exists an
$\epsilon_N>0$ so that for all $0<\epsilon<\epsilon_N$, there
exist at least two stationary solutions which intersect $u=0$
exactly $N$ times. (The cited works impose Neumann boundary
conditions: it is a simple generalization to periodic boundary
conditions.)

Via Theorem~\ref{thm_infinite}, we have that for any such $h$ and
any small $\epsilon$, this equation admits an infinite collection
of stationary periodic curves; furthermore, there is a lower bound
of $N$ on the number of $1$-periodic solutions.
}\end{example}

\begin{example}
\label{ex_Stat}{\em
Consider the following equation

\begin{equation}
   u_t = u_{xx} - \frac{5}{8}\sin 2x\, u_x
   + \frac{\cos x}{\cos x + \frac{3}{\sqrt{5}}}
   u(u^2-1),
\end{equation}
with $x\in S^1=\re/2\pi \zet$.

Eqn.~(\ref{eq_Stat}) is a weighted exact system with Lagrangian
\begin{equation}
     L = e^{-\frac{5}{16}\cos2x}\left(
    \frac{1}{2}u_x^2-\frac{\cos x}{\cos x + \frac{3}{\sqrt{5}}}
    \frac{(u^2-1)^2}{4}\right) ,
\end{equation}
where by ``weighted exact'' we mean (cf. \rmref{wexact})
\begin{equation}
    u_t = e^{\frac{5}{16}\cos2x}\left[
        \frac{d}{dx}\frac{\partial L}{\partial u_x}
        -\frac{\partial L}{\partial u}\right]
\end{equation}

One checks easily that there are stationary solutions $u=\pm 1$
and $u_\pm=\pm\frac{1}{2}\left(\sqrt{5}\cos x +1\right)$, as in
Fig.~\ref{fig_Stat}. These curves comprise a skeleton
$\v=\{-1,u_-, u_+,+1\}$ which can be discretized to yield the
skeleton of Example~\ref{ex_indexex}. From the computation of the
index there, this skeleton forces a stationary solution of the
braid class indicated in Fig.~\ref{fig_indexex}[left]: of course,
this is detecting the obvious stationary solution $u=0$.

What is more interesting is the fact that one can take periodic
extensions of the skeleton and add free strands in a manner which
makes the relative braid spatially non-periodic. Let us describe a
family of proper and bounded relative braid classes. Let $\v^n$ be
the $n$-fold periodic extension of $\v$ on $[0,n-1]$ and consider
a single free strand $-1<u(x)<1$ that links with $\v^n$ as
follows: on each interval $[k,k+1]$, $k=0...n-1$, we choose one of
three possibilities:
\begin{enumerate}
\item[(a)] $\iota(u,u_-)=0$ and $\iota(u,u_+)=2$,
\item[(b)] $\iota(u,u_-)=2$ and $\iota(u,u_+)=2$, or
\item[(c)] $\iota(u,u_-)=2$ and $\iota(u,u_+)=0$.
\end{enumerate}
Define a symbol sequence $\sigma = (\sigma_1...\sigma_n)$, where
$\sigma_i \in \{a,b,c\}$. Every symbol sequence except for
$\sigma=(a...a)$ and $\sigma=(c...c)$, defines a proper and
relative braid class $\{\u_\sigma\rel\v\}$.

To compute the invariant, we discretize. Choose the discretization
$d=2n$ on $[0,n]$.\footnote{This is admissible for the skeleton
$\v^n$ and is large enough to yield the correct index
computation.} Fig.~\ref{fig_indexex}[right] shows an example. In
\S\ref{sec_inv} the index was computed:
$$
\HH(\u_\sigma\rel\v)=h(\disc_{2n}\u_\sigma\rel \disc_{2n}\v^n)\simeq S^k ,
$$
where $k = \#\{b\in\sigma\}$. Therefore, $P_\tau(\HH) = \tau^k$.

The Morse inequalities now imply  that for each $n>0$ there exist
at least $3^n -2$ different stationary solutions. This information
can be used again to prove that the time-$2\pi$ map of the
stationary equation has positive entropy.

\begin{figure}[hbt]
\begin{center}
\psfragscanon
\psfrag{p}[r][]{$+1$}
\psfrag{m}[r][]{$-1$}
\includegraphics[width=2in]{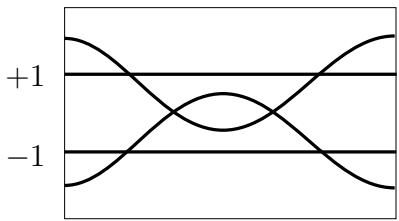}
\caption{The skeleton of stationary solutions for
Eqn.~(\ref{eq_Stat}) forces an infinite collection
of additional solutions which grows exponentially
in the number of strands employed.}
\label{fig_Stat}
\end{center}
\end{figure}
}\end{example}

\begin{example}{\em
\label{ex_InfMany} The following class of examples is very general
and includes Example \ref{ex_Matano} as a special case. One says
that \rmref{paracont} is {\em dissipative} if
\begin{equation}
    u\,f(x,u,0,0) \to -\infty {\mbox{ as }} \abs{u}\to+\infty
\end{equation}
uniformly in  $x \in S^1$.

\begin{theorem}
\label{thm_infinite} Let $f$ be dissipative and satisfy {\bf
(f1)}-{\bf (f2)}. If $\v$ is a nontrivially braided stationary
skeleton (i.e., $\iota(\v)\neq 0$), then there are infinitely many
braid classes represented as stationary solutions to
\rmref{paracont}. Moreover, the number of braid classes for which
$\u$ consists of just one strand is bounded from below by $2\lceil
{\iota/2}\rceil$, where $\iota$ is the maximal number of
intersections between two strands of $\v$.
\end{theorem}
\proof Given the assumptions one can find $c_+>0$ and $c_-<0$ such
that $\pm f(x,c_\pm,0,0)>0$, and
$$
c_- < v^\alpha(x) < c_+,
$$
for all strands $v^\alpha$ in $\v$. Using discrete enclosure via
sub/super solutions, Lemma \ref{help} in Appendix C yields
solutions $u_+$ and $u_-$ such that
$$
c_- < u_-(x) < v^\alpha (x) < u_+(x) < c_+,
$$
for all $\alpha$. Assume without loss of generality that all
strands in $\v$ are 1-periodic (if not, one can take an
appropriate covering of $\v$). For the sake of convenience we may
assume that $x\in S^1\equiv \re/\zet$. Select two intersecting
strands which form the braid $\w=\{v^{\alpha_1},v^{\alpha_2}\}$,
and set $\iota(\w) = \# \{ \hbox{\rm intersections between}
~v^{\alpha_1} ~\hbox{\rm and}~ v^{\alpha_2} \}$. Consider the
skeleton $\z =\{u_-,v^{\alpha_1},v^{\alpha_2},u_+\}$ and a free
strand $u(x)$ --- with $u(x+1)=u(x)$ --- that links with $\z$ as
follows: (1) $u_-(x) \le   u(x) \le u_+(x)$, (2) for some $k>0$,
$\iota(u,v^{\alpha_1})=\iota(u,v^{\alpha_2})=2k <\iota(\w)$. These
two hypotheses describe the relative braid class $\{\u\rel \z\}$,
which clearly is a proper and bounded class and therefore has a
well-defined homotopy braid index $\HH$. The index $\HH$ is an
invariant of the braid class and it can be computed for instance
by studying a specific system of which all solutions are known.

Consider the equation $\epsilon^2 u_{xx} + u -u^3=0$. If we choose
$\epsilon = \left(\pi(\iota(\w)+1)\right)^{-1}$, then there exists
a periodic solution $v_1(x)$ with period $T = 2/\iota(\w)$. Define
$v_2(x)=v_1\bigl(x- (1/\iota(\w))\bigr)$; then if we consider
$v_1$ and $v_2$ on the interval $[0,1]$, it follows that
$\iota(v_1,v_2) = \iota(\w)$. The skeleton $\z'=\{-1,v_1,v_2,+1\}$
is now topologically equivalent to $\w$. Moreover, the equation
$\epsilon^2 u_{xx} + u -u^3=0$ has a unique solution $u$ which has
the right linking properties with the skeleton $\z'$: $-1<u(x)<1$,
and $\iota(u,v_1)=\iota(u,v_2) = \iota(u,v^{\alpha_1})
=\iota(u,v^{\alpha_2})<\iota(\w)$. Therefore, $\{\u\rel \z'\}$ and
$ \{\u\rel\z\}$ are topologically equivalent. As in
\cite{GVV,Ang3} the invariant set $\inv(\{\u\rel \z'\})$ of the
equation
$$
u_t =\epsilon^2 u_{xx} + u -u^3,
$$
is given by $\inv(\{\u\rel \z'\}) =\{u(x+\phi)~|~\phi\in\re\}$,
which represents a hyperbolic circle of stationary strands. Its
unstable manifold has dimension $\iota(u,v_1)=2k$ and therefore
its Morse polynomial is given by $\tau^{2k-1} (1+\tau)$. Since
this captures the entire invariant, it follows that $P_\tau(\HH) =
\tau^{2k-1} (1+\tau)$: see also \cite{Ang3} for details.

From the invariant $\HH$ and Theorem \ref{thm_Main} we deduce that
if \rmref{paracont} is dissipative and exact it has at least
$\lceil {\iota(\w)\over 2}\rceil$ pairs of $1$-periodic solutions.
One finds infinitely many stationary braids by allowing periods
$2n$. Indeed, take the periodic extension $\w^{2n}$. Then for any
$k$ satisfying $2k < \iota(\w^{2n}) = 2n\iota(\w)$ we find a
$2n$-periodic solution. By projecting this to the interval $[0,1]$
we obtain a multi-strand stationary braid for \rmref{paracont}. As
a matter of fact for each pair $p,q$, with $q<p$ and ${\rm
gcd}(p,q)=1$, there exists at least two distinct periodic
solutions $u^1_{p,q}$ and $u^2_{p,q}$, by setting $k=\iota(\w)q$
and $n=p$. This infinity of solutions enshrouds the set $\queu
\cap (0,1)$. \qed }\end{example}

\section{Examples: time periodic solutions}
\label{sec_exper}

This is a longer example of a very general forcing result for
time-periodic solutions.

\begin{example}
\label{ex_Pers}{\em
Consider equations of the following type
    \begin{equation}
    \label{eq_Pers}
        u_t = u_{xx}+ug(u)+u_xh(x,u,u_x),\quad x\in \real/\zet,
    \end{equation}
where the non-linearity is assumed to satisfy {\bf (f2)}, i.e. $h$
has sub-linear growth in $u_x$ at infinity. Moreover, assume that
$g$ and $h$ satisfy the hypotheses:
\begin{itemize}
\item[\bf{(g1)}] $g(0) >0$, and $g$ has at least one
    positive and one negative root;
\item[\bf{(g2)}] $h>0$ on $\{ u u_x\neq 0 \}$.
\end{itemize}
\begin{theorem}
Under the hypotheses above \rmref{eq_Pers} possesses an infinite
collection of time-periodic solutions all with different braid
classes.
\end{theorem}
}\end{example}

\begin{proof}
Consider first the perturbed equation,
\begin{equation}
\label{peqn}
    u_t = u_{xx}+ug(u)+\alpha_\epsilon u_xh(x,u,u_x),
\end{equation}
where $\alpha_\epsilon=0$ for $\sqrt{u^2+u_x^2}\in[0,\epsilon]$
and $\alpha_\epsilon=1$ for $\sqrt{u^2+u_x^2}\ge 2\epsilon$. For
$\epsilon>0$ \rmref{peqn} has small stationary solutions
$u_\epsilon(x)$ which oscillate about $u=0$. We can choose this
$u_\epsilon$ and a relatively prime pair of integers $p,q \in
\nat$ such that $u_\epsilon(x+p)=u_\epsilon(x)$ and
$\sqrt{g(0)}/2\pi\leq q/p$ is arbitrarily close to $q/p$. The
integer $q$ represents the number of times the oscillation fits on
the interval $[0,p]$.

We use {\bf (g1)} to build a skeleton for \rmref{peqn}. Let $a_+$
and $a_-$ denote positive and negative roots of $g$, and consider
the skeleton $\v=\{v^1,v^2,v^3,v^4\}$ on $\real/p\zet$ with
$v^1(x)=a_-$, $v^2(x)=a_+$, $v^3(x) = u_\epsilon(x)$, and $v^4(x)
= u_\epsilon(x-p/2q)$. Clearly $\iota(v^3,v^4) = 2q$. Define the
relative braid class $\{\u\rel\v\}$ as follows; $\u=\{u\}$ is a
(1-strand) braid satisfying $a_-<u(x)<a_+$ and
$\iota(u,v^3)=\iota(u,v^4)=2r < 2q$. This braid class is proper
and bounded, and its homotopy invariant $\HH$ was computed in the
previous section:
$$
    P_\tau(\HH) = \tau^{2q-1}(1+\tau).
$$

We claim that for $0<\epsilon\ll 1$ there are no stationary
solutions in $\{\u\rel\v\}$. Suppose that $u$ is stationary. One
checks that the function
$$
    H(u,u_x) := {1\over 2} u_x^2
            + u\int_0^ug(s)ds
            -\int_0^u\int_0^s g(r)dr\,ds
$$
has derivative
$$
    \frac{d}{dx}H = -\alpha_\epsilon u_x^2 h(x,u,u_x).
$$
This term is nonpositive by {\bf(g2)} and not identically zero by
the fact that $u$ cannot be close to a constant (thanks to the
intersection numbers). The periodic boundary condition leads to
the desired contradiction.

Since $\u$ is a 1-strand braid it follows from Theorem
\ref{thm_Main2} that $\{\u\rel\v\}$ contains a $t$-periodic
solution to \rmref{peqn}. By lifting the equation to the interval
$[0,kp]$, $k\in \nat$, we obtain different periodic solutions for
each $r<kq$, which shows that there are $t$-periodic solutions for
infinitely many different braid classes: see \cite[Lem. 43]{GVV}
for details. What remains is to show that these periodic solutions
to \rmref{peqn} persist in the limit $\epsilon \to 0$. We need to
show that the limits obtained are not equal to the zero solution.
We use an argument similar to that of Angenent \cite{Ang3}.

Linearize \rmref{peqn} around $u=0$. This leads to the linear
operator $L={d^2/dx^2}+g(0)$ on $L^2(\real/p\zet)$. The spectrum
of $L$ is given by the eigenvalues $\lambda_n = -4\pi^2 n^2/p^2 +
g(0)$, for $n=0,1,...$. Since $\sqrt{g(0)}/2\pi\le q/p$ it holds
that $\lambda_n>0$ for all $n<q$, and $\lambda_n \le 0$ for $n\ge
q$. This yields a (spectral decomposition) splitting of
$L=L_++L_-$. The evolution on the set
$I=\{\psi~|~\iota(\psi,0)=2r<2q\}$ is then dominated by the linear
operator $L$ for $\Vert\psi\Vert_{L^2}$ small. Therefore, the
function $B(\psi) = {1\over 2}(\psi,L_+\psi)_{L^2}$ satisfies
\begin{equation}
\label{eq_B}
\frac{d}{dt} B(\psi) = (L_+\psi,L_+\psi)_{L^2}
    + o(\Vert\psi\Vert_{L^2}) . 
\end{equation}
for all $\psi \in I$. For $u_\epsilon(t,x)$ a periodic solution
with sufficiently small $L^2$ norm, \rmref{eq_B} implies that
$\frac{\partial}{\partial t}B>0$, a contradiction of periodicity.
Thus we conclude that the $\epsilon\to 0$ limits do not collapse
to zero.
\end{proof}
\begin{remark}{\em
\label{rem_General} The form of Eqn.~(\ref{eq_Pers}) is not the
most general form possible. Certainly, having $h$ strictly
negative on $\{u u_x\neq 0\}$ is also permissible. With work, the
diligent reader may verify that allowing the $u_{xx}$ term to vary
as per {\bf (f1)} does not change the nature of the results.
}\end{remark}

\section{Proofs: Forcing stationary solutions}
\label{sec_proof}

\subsection{Discretization of the equation}
From hypothesis {\bf (f1)} we obtain an estimate for $f$ of the
form
\begin{equation}
\label{eq_est}
    f(x,u,v,0) + a_-(w) w \le f(x,u,v,w) \le f(x,u,v,0) + a_+(w) w,
\end{equation}
for all $x \in [0,1]$, and $u,v,w \in \re$,
where $a_-(s) = \lambda^{-1}$ for $s\le0$, $a_-(s)=\lambda$ for $s\ge 0$, and
$a_+(s) = \lambda $ for $s\le0$, $a_+(s)=\lambda^{-1}$ for $s\ge 0$.

Consider a braid $\u$ of $n$ strands. For the remainder of this
section, we work with individual strands $u=u^\alpha$, suppressing
the superscripts for notational aesthetics.

We discretize \rmref{paracont} in the standard manner. Choose a
step size $1/d$, for $d \in \nat$, and define $u_i:=u(i/d)$. We
approximate the first derivative $u_x(i/d)$ by $\d u_i :=
d(u_{i+1}-u_i)$ and the second derivative $u_{xx}(i/d)$ by $\dd
u_i :=d^2 (u_{i+1}-2u_i+u_{i-1})$.
\begin{lemma}\label{back}
Let $\u$ be a stationary braid for \rmref{paracont}, then
\begin{equation}
      \epsilon_i(d) :=
      f\Bigl({i\over d},u_i,\Delta u_i,\Delta^2u_i)\Bigr) \longrightarrow 0,
\end{equation}
as $d\to \infty$ uniformly in $i$.
In particular, $|\epsilon_i(d)|\le C/d$.
\end{lemma}
\begin{proof}
From Appendix A it follows that each strand $u$ of a stationary
solution to \rmref{paracont} is $C^3$. A Taylor expansion yields
\begin{eqnarray*}
    \d u_i -u_x &=&
    d\cdot \bigl(u((i+1)/d)-u_i(i/d)\bigr)-u_x(i/d)\\
    &=& d\cdot  R^2_{1/d}(i/d) = {1\over 2} u_{xx}(y)/d,
    \,\, {\mbox{ for some }} y \\
    \dd u_i -u_{xx}&=&
    d^2\cdot \bigl(u((i+1)/d)-2u_i(i/d)
        +u((i-1)/d)\bigr) -u_{xx}(i/d) \\
    &=& d^2 \cdot R^3_{1/d}(i/d)
    = {1\over 6} u_{xxx}(y)/d.
    \,\, {\mbox{ for some }} y
\end{eqnarray*}
For $x_0 := i/d $ it therefore holds that
\[
\bigl| \d u_{d\cdot x_0} - u_x(x_0)\bigr| \le C/d, \quad\quad
\bigl| \dd u_{d\cdot x_0} - u_{xx}(x_0)\bigr| \le C/d,
\]
with $C$ independent of $x_0$. Since $f$ is $C^1$ the desired
result follows. A more detailed asymptotic expansion for
$\epsilon_i$ is obtained as follows:
\begin{eqnarray*}
\epsilon_i(d) &=& f(i/d,u_i, \d u_i,
 \dd u_i) - f(i/d,u(i/d),u_x(i/d),u_{xx}(i/d))\\
 &=&\partial_{u_x} f(i/d,u(i/d),u_x(i/d),u_{xx}(i/d))
 \bigl[ \d u_i - u_x(i/d)\bigr]\\
 &&+\partial_{u_{xx}} f(i/d,u(i/d),u_x(i/d),u_{xx}(i/d))
 \bigl[ \dd u_i - u_{xx}(i/d)\bigr] + R_2 .
\end{eqnarray*}
 From the weak form of Taylor's Theorem the remainder term $R_2$ satisfies
$|R_2| = o(1/d)$. Combining this with the estimates obtained above
we derive that $|\epsilon_i(d)| \le C/d$, thus completing the
proof.
\end{proof}

The next step is to ensure that the $d$-point discretization of
$\u$ is in fact a solution of an appropriate parabolic recurrence
relation.
\begin{lemma}
\label{approx} Let $\u \in \Conf^n$, and let $d$ be an admissible
discretization for $\u$. Then for any sequence
$\{\epsilon_i^\alpha\}$ with $i=0,..,d$ and $\alpha=1,..,n$, there
exists a parabolic recurrence relation $\EE_i^d$ satisfying
\[
    \EE^d_i(u_{i-1}^\alpha,u_i^\alpha,u_{i+1}^\alpha)
    = -\epsilon_i^\alpha,
\]
where $u_i^\alpha = u^\alpha(i/d)$. In addition, $|\EE^d_i(r,s,t)|
\le C\max_{i,\alpha} |\epsilon^\alpha_i|$ for all $|r|, |s|, |t|
\le 2 \max_{i,\alpha} |u_i^\alpha|$ and some uniform constant $C$
depending only on $\u$.
\end{lemma}
\begin{proof}
This proof is a straightforward extension of \cite[Lemma 55]{GVV},
in which $\epsilon_i^\alpha\equiv 0$.
\end{proof}

Let $\v\in \Conf^m$ be stationary for \rmref{paracont} and let
$\{\u\rel\v\}$ be a bounded proper braid class with $d$ a
sufficiently large discretization period. We now construct a
parabolic recurrence relation for which the discrete skeleton
$\disc_d\v$ is stationary. Combining the Lemmas~\ref{back} and
\ref{approx}, the recurrence relation defined by
\begin{equation}
    \R_i^d(u_{i-1},u_i,u_{i+1}) :=
    f(i/d,u_i,\d u_i, \dd u_i) + \EE_i^d(u_{i-1},u_i,u_{i+1}),
\end{equation}
has $\disc_d\v$ as a stationary solution. The above construction
works for any $d'\ge d$. For a given $d$ the recurrence relation
$\R_i^d$ is considered on the compact set ${\rm cl}[\u\rel\v]$,
which implies that $|u_i| < 2\max_{i,\alpha} |v_i|$. To verify the
parabolicity of $\R_i^d$, we compute the derivatives.
 From hypothesis {\bf (f1)} and the parabolicity of $\EE_i^d$, we obtain:
\begin{equation}
    \partial_1\R^d_i = {\partial  f \over \partial u_{xx}} \cdot  d^2
    + \partial_1\EE^d_i \ge \lambda \cdot d^2 >0.
\end{equation}
Furthermore,
\[
    \partial_3\R^d_i ={\partial  f \over \partial u_{xx}} \cdot  d^2 +
    {\partial  f \over \partial u_{x}} \cdot d
    + \partial_3\EE^d_i,
\]
which does not yet prove parabolicity since no estimates for
$\partial_{u_x} f$ are given.
However,  utilizing hypothesis {\bf (f2)}, we have that
\begin{eqnarray*}
\R_i^d(u_{i-1},u_i,u_{i+1}) &=&
    f(i/d,u_i,\d u_i, \dd u_i) + \EE_i^d(u_{i-1},u_i,u_{i+1})\\
&\ge& f(i/d,u_i,\d u_i,0) + a_-(\dd u_i) \dd u_i + \EE_i^d\\
&\ge& -C -C|\d u_i|^\gamma +\lambda \dd u_i\\
&\ge& -C -Cd^\gamma  +\lambda d^2(u_{i+1} -2u_i +u_{i-1}),
\end{eqnarray*}
which shows that $\R_i^d$ is an increasing function of $u_{i+1}$,
provided that $d$ is large enough. This relies on the fact that
the braid class is bounded. Periodicity of $f$ implies that
$R^d_{i+d} = \R^d_i$. Initially $\R_i^d(r,s,t)$ is defined for
$|r|, |s|, |t| \le 2\max_{i,\alpha} |v_i|$. It is clear that
$\R_i^d$ can easily be extended to a parabolic recurrence relation
on all of $\re^3$.

\subsection{Convergence to a stationary solution}

Choose $d_*$ large enough such that $\R_i^d$ is parabolic for all $d\ge d_*$.
Let $\{u_i^{\alpha,d}\}$ be a sequence of braids
which are solutions of
\begin{equation}
\label{mainrec}
    \R_i^d(u_{i-1},u_i,u_{i+1})
    = f(i/d,u_i,\d u_i, \dd u_i) +
    \EE_i^d(u_{i-1},u_i,u_{i+1}) =0,
\end{equation}
and which satisfy the uniform estimate $|u_i^{\alpha,d}| \le C$
for all $d$. For notational simplicity, we omit the discretization
period $d$ and write $u_i$ instead of $u_i^{\alpha,d}$ in what
follows. The discretization index will be clear from the range of
the index $i$.

\begin{lemma}\label{estI}
Let $\{u_i\}$ satisfy
$\R_i^d(u_{i-1},u_i,u_{i+1}) =0$
and $|u_i|\le C$ as $d\to \infty$. Then
\begin{equation}
    \sum_{i=0}^d {1\over d} |u_i|^2 \le C
        \quad ; \quad
    \sum_{i=0}^d d\cdot |u_{i+1} -u_i|^2 \le C,
\end{equation}
with $C$ independent of $d$.
\end{lemma}
\begin{proof}
For each strand $\alpha$ it holds that either $u_{i+d} = u_i$, or
$u_{i+kd} = u_i$ for some $k$. Since there are only finitely many
strands, the constant $k$ can be chosen uniformly for all
$\alpha$. Therefore we assume without loss of generality that the
first equality holds. The first estimate immediately follows from
the uniform bound on $u_i$.

From \rmref{mainrec} it follows that
\[
    f(i/d,u_i,\d u_i, \dd u_i) =-
    \EE_i^d(u_{i-1},u_i,u_{i+1}).
\]
Multiply the above equation by $u_i/d$, then from \rmref{eq_est}
it follows that
\[
    -\dd u_i \cdot {1\over d} u_i \le
    {1\over a_+} f(i/d,u_i,\d u_i,0)
    {1\over d} u_i + {1\over a_+}
    \EE_i^d(u_{i-1},u_i,u_{i+1})
    {1\over d} u_i,
\]
for $u_i>0$, and
\[
    -\dd u_i \cdot {1\over d} u_i \le +
    {1\over a_-} f(i/d,u_i,\d u_i,0)
    {1\over d} u_i +{1\over a_-}
    \EE_i^d(u_{i-1},u_i,u_{i+1})
    {1\over d} u_i,
\]
for $u_i<0$. From the periodic boundary conditions it follows that
\[
    -\sum_{i=0}^d
    \dd u_i \cdot {1\over d} u_i
    = \sum_{i=0}^d d\cdot|u_{i+1} -u_i|^2.
\]
Combining the above estimates and using {\bf (f2)} we obtain
\begin{eqnarray*}
    \sum_{i=0}^d d\cdot|u_{i+1} -u_i|^2
&\le&
    \sum_{i=0}^d \Bigl( {\lambda\over d}
    |f(i/d,u_i,\d u_i,0)| |u_i|
    + {\lambda\over d} |\EE_i^d||u_i|\Bigr)
\\
&\le&
    C/d +\sum_{i=0}^d {\lambda\over d} \Bigl[ C_\epsilon
    + \epsilon d^2|u_{i+1} -u_i|^2\Bigr]|u_i|
\\
&\le&
    C_\epsilon + \epsilon C
    \sum_{i=0}^d d\cdot |u_{i+1} -u_i|^2,
\quad\quad
\hbox{for any}~~\epsilon>0.
\end{eqnarray*}
Choosing $\epsilon$ small enough yields the second estimate.
\end{proof}

Define $\phi_d := \pl(\{u_i\})$. Then
$\| \phi_d\|^2_{L^2} \le \sum_{i=0}^d {1\over d} |u_i|^2$,
and $\| {d\over dx} \phi_d\|^2_{L^2} =
 \sum_{i=0}^d d\cdot|u_{i+1} -u_i|^2$.
Due to the  uniform estimates we obtain the Sobolev bound
$\|\phi_d \|_{H^{1,2}} \le C$, with $C$ independent of $d$.
Therefore, $\phi_{d_n}$ converges to some function $u \in
C^0([0,1])$, as $ d_n \to \infty$.

\begin{lemma}\label{estII}
Let $\{u_i\}$ satisfy
$\R_i^d(u_{i-1},u_i,u_{i+1}) =0$ and
$|u_i|\le C$ as $d \to \infty$. Then
\begin{equation}
  \sum_{i=0}^d d\cdot |\d u_{i} -\d u_{i-1}|^{2/\gamma}
    \le C,
\end{equation}
with $C$ independent of $d$.
\end{lemma}
\begin{proof}
As in the proof of Lemma~\ref{estI}, we have
\[
    \dd u_i \le -
    {1\over a_-} f(i/d,u_i,\d u_i,0)
    -{1\over a_-} \EE_i^d(u_{i-1},u_i,u_{i+1}),
\]
for $ \dd u_i >0$, and
\[
    -\dd u_i \le
    {1\over a_+} f(i/d,u_i,\d u_i,0)
    +{1\over a_+} \EE_i^d(u_{i-1},u_i,u_{i+1}),
\]
for $\dd u_i <0$. Combining these estimates with {\bf (f2)} we
obtain
\begin{equation}
\label{inter}
    d\cdot | \d u_i - \d u_{i-1} |
    \le C + C |\d u_i|^\gamma.
\end{equation}
Therefore
\begin{eqnarray*}
    \sum_{i=0}^d d^{{2\over \gamma} -1} |
    \d u_i - \d u_{i-1} |^{2/\gamma}
    &\le&
    C+C \sum_{i=0}^d {1\over d} |\d u_i|^2
\\
    &\le&
    C \quad \hbox{by Lemma \ref{estI},}
\end{eqnarray*}
which is the desired estimate.
\end{proof}
Set $\psi_d := \pl(\{\d u_i\})$. Then $\| \psi_d \|^2_{L^2} \le
\sum_{i=0}^d {1\over d} |\d u_i|^2 \le C$, and $\| {d\over dx}
\psi_d\|^2_{L^{2/\gamma}} =
 \sum_{i=0}^d d^{{2\over \gamma} -1} \cdot|\d u_{i+1}
 -\d u_i|^{2/\gamma} \le C$.
This implies that $\|\psi_d\|_{H^{1,2/\gamma}} \le C$,
independent of $d$. Therefore there exists a subsequence
$\psi_{d_n}$ converging to some function $v \in C^0([0,1])$.

\begin{lemma}\label{estIII}
Let $\{u_i\}$ satisfy $\R_i^d(u_{i-1},u_i,u_{i+1}) =0$, and
$|u_i|\le C$ as $d \to \infty$, then
\[
    |\d u_i| \le C,
        \quad\quad
    |\dd u_i| \le C,
\]
with $C$ independent of $d$.
\end{lemma}
\begin{proof}
The first estimate follows from the fact that $\|\psi_d\|_{C^0}
\le C$, hence $|\d u_i| \le C$. For the second estimate we use
\rmref{inter}. The uniform bound on $\d u_i$ then yields a uniform
bound on $\dd u_i$.
\end{proof}
Finally, we require an estimate on
$\Delta^3 u_i = d\cdot (\dd u_{i+1} - \dd u_i)$.

\begin{lemma}\label{estIV}
Let $\{u_i\}$ satisfy $\R_i^d(u_{i-1},u_i,u_{i+1}) =0$ and
$|u_i|\le C$ as $d \to \infty$. Then
\begin{equation}
    |\Delta^3 u_i| = d\cdot |\dd u_{i+1}
    - \dd u_i| \le C,
\end{equation}
with $C$ independent of $d$.
\end{lemma}
\begin{proof}
Since $\R_{d+1}^d -\R_d^i=0$ it follows from the definition of
$\R_i^d$ that
\begin{eqnarray*}
&~&
    f((i+1)/d,u_{i+1},\d u_{i+1},\dd u_{i+1}) -
    f(i/d,u_i,\d u_i, \dd u_i)
\\
&=&
    -\EE_{i+1}^d(u_{i},u_{i+1},u_{i+2})
    + \EE_{i}^d(u_{i-1},u_i,u_{i+1}).
\end{eqnarray*}
Using Taylor's theorem we obtain that
\begin{eqnarray*}
    \partial_x f(i/d,u_i,\d u_i,\dd u_i){1\over d}
    &+&
    \partial_u f(i/d,u_i,\d u_i, \dd u_i)
    (u_{i+1}-u_i)
\\
    &+&
    \partial_{u_x} f(i/d,u_i,\d u_i, \dd u_i)
    (\d u_{i+1} - \d u_i)
\\
    &+&
    \partial_{u_{xx}} f(i/d,u_i,\d u_i, \dd u_i)
    (\dd u_{i+1} - \dd u_i)
\\
        &=&
    - (\EE_{i+1}^d -\EE_i^d)
    - R_2(i/d,u_i,\d u_i, \dd u_i).
\end{eqnarray*}
For $\Delta^3 u_i$ this implies
\begin{eqnarray*}
    \partial_{u_{xx}} f(i/d,u_i,\d u_i, \dd u_i)
    \Delta^3 u_i
    &=&
    -\partial_x f(i/d,u_i,\d u_i, \dd u_i)
\\
    &-&
    \partial_u
    f(i/d,u_i,\d u_i, \dd u_i)\d u_i
\\
    &-&
    \partial_{u_x}
    f(i/d,u_i,\d u_i, \dd u_i)\dd u_i
\\
    &-&
    \d \EE_i^d -R_2(i/d,u_i,\d u_i, \dd u_i)d.
\end{eqnarray*}
By Lemma \ref{estIII} the right hand side is uniformly bounded in
$d$. Using {\bf (f1)} then yields the desired estimate on
$\Delta^3 u_i$.
\end{proof}
Define $\chi_d := \pl(\{\dd u_i\})$. From Lemma \ref{estIV} we
then derive that $\| {d\over dx} \chi_d\|_{L^\infty} \le C$.
Therefore $\chi_{d_n}$ converges to some limit function $w$ in
$C^0([0,1])$.

From Lemmas \ref{estI}, \ref{estII}, and \ref{estIV} it follows
that the functions $\phi_d$, $\psi_d$ and $\chi_d$ converge to
function $u$, $v$ and $w$ respectively, with the anchor points
being solutions of $\R_i^d =0$.

The following lemma relates discretized braids to stationary
braids in $\{ \u\rel \v\}$.
\begin{lemma}\label{conv}
Let $u, v, w \in C^0([0,1])$ and let $\{u_i^d\}_{i=0}^d$ be
sequences whose PL interpolations satisfy
\[
    \pl (u_i^d)  \to u,\quad
    \pl (\d u_i^d)  \to v, \quad
    \pl (\dd u_i^d)  \to w,
\]
in $C^0([0,1])$ as $d\to \infty$. If $\R_i^d(u_{i-1}^d,
u_i^d,u_{i+1}^d) = 0$ and $|\Delta^3 u_i| \le C$, then $u \in
C^2([0,1])$, $u_x=v$, and $u_{xx}=w$ satisfying
$f(x,u,u_x,u_{xx})=0$ pointwise on $[0,1]$.
\end{lemma}
\begin{proof}
We start with the estimate $|\phi'_d - \psi_d| \le {1\over d}
\max_{0\le i\le d}|\dd u_i|\to 0$ uniformly as $d \to \infty$.
This implies that $ \psi_d \to v$ in $C^0([0,1])$. The same
estimate holds for $|\psi'_d - \chi_d| \le {1\over d} \max_{0\le
i\le d}| \Delta^3 u_i| \to 0$ uniformly as $d \to \infty$. Hence
we deduce that $\chi_d \to w$ in $C^0([0,1])$. From the definition
of derivatives it now follows that $\|D_{1/d} u -v\|_{L^\infty}
\to 0$, and $\|D_{1/d} v -w\|_{L^\infty} \to 0$; thus $v=u_x$ and
$w=u_{xx}$. From Lemma~\ref{back} we deduce that
$f(x,u,u_x,u_{xx})=0$.
\end{proof}

Note that Lemma~\ref{lem_A} in Appendix A implies further that
$u\in C^3([0,1])$.

\subsection{Proof of Theorem~\ref{thm_Main}}\label{sec_Main}

Given $P_{-1}\HH\neq 0$, the existence of a single stationary
solution is argued as follows. Choose $d_*$ large enough. Then
from Theorem \ref{thm_GVV2} it follows that \rmref{mainrec} has a
discrete braid solution $\{u_i^{\alpha,d}\}$ for all $d\ge d_*$.
The boundedness of the braid class implies that the sequence
$\{u_i^{\alpha,d}\}$ satisfies
$$
    |u_i^{\alpha,d}|\le C,\quad\quad\forall~
    i,\alpha,~{\rm and}~\forall~d\ge d_*.
$$
Lemmas~\ref{back}-\ref{conv} imply that as $d\to\infty$ one
obtains a stationary braid $\u=\{u^\alpha\}$ whose strands
$u^\alpha$ satisfy the equation
$f(x,u^\alpha,u_x^\alpha,u_{xx}^\alpha)=0$. Since the skeletal
strands $\disc_d(v^\beta)$ converge to $v^\beta$ by construction
and the pairwise intersection numbers are the same for all $d$, we
have in the limit a solution to \rmref{paracont} in the correct
braid class.

It remains to determine multiplicity in the case of {\bf (f3)}.
The difficulty lies in dealing with degenerate critical points:
one proceeds using the standard tools of critical groups and
Gromoll-Meyer pairs. We refer the interested reader to
\cite{Chang} for detailed definitions. For the remainder of the
proof, we will characterize Morse data of critical points $\u$ via
the Poincar\'e polynomial $P_\tau(\u)$. For a nondegenerate
critical point, this is a polynomial of the form
$P_\tau(\u)=\tau^{\mu(\u)}$, where $\mu$ is the Morse index. For
degenerate critical points, $P_\tau$ is defined via certain
homology groups \cite{Chang}.

In the gradient case one has the action $\A$ on the space
$\overline{\Conf^n}$ defined as follows:
\[
    \A(\u) = \sum_{\alpha=1}^n \int_0^1 L(x,u^\alpha,u_x^\alpha)dx,
    \quad \u \in \overline{\Conf^n},
\]
and the discretized action on $\overline{\DConf_d^n}$ defined by
\[
    \A_d(\u) = \sum_{\alpha=1}^n\sum_{i=0}^d
    L\Bigl({i\over d},u_i^\alpha,\d u_i^\alpha\Bigr) +
        \sum_{\alpha=1}^n\sum_{i=0}^d a_i(u_i^\alpha,\d u_i^\alpha),
    \quad \u \in \overline{\DConf^n_d},
\]
where the $a_i$ are small perturbations guaranteeing that $\disc_d
\v$ is a critical skeleton for each $d\ge d_*$.\footnote{We omit
the superscript $d$ in the notation for $u_i^\alpha$.} These can
be constructed as in Lemmas~\ref{back}-\ref{approx} so as to
satisfy the same estimates. It follows immediately from
\rmref{wexact} that $\R_i^d = -\partial_{u_i} \A_d$.

Assume without loss of generality that $\A$ has finitely many
critical points $\u_i\rel\v$ so that all critical points are
isolated. We have shown earlier in this section that as
$d\to\infty$, critical points of $\A_d$ converge to a critical
point of $\A$. We will factor this convergence through a sequence
of nondegenerate Morse functionals in order to extract forcing
data.

One may perturb $\A$ on a neighborhood of the critical points to
$\A^\epsilon$ which is Morse on the braid class $\{\u\rel\v\}$.
Next, discretize $\A^\epsilon$ to yield functionals
$\A^\epsilon_d$. Our convergence results imply that
$\A^\epsilon_d$ is Morse for $d$ sufficiently large. Indeed, if
$\{\u^d\}$ is a sequence of critical points of $\A^\epsilon_d$,
then $\pl(\u^d)$ converges in $\overline{\Conf^n}$ to a critical
point of $\A$. The same holds for the eigenfunctions and
eigenvalues of the linearized functional, which implies that
$\A^\epsilon_d$ is Morse for $d$ large enough. Uniform estimates
on the remainder terms of $\A^\epsilon$ and $\A_d^\epsilon$ then
yield uniformity in the distance between the critical points of
$\A^\epsilon_d$ for all $d$ large. To be more precise
$\dist_{\overline\DConf_d^n}(\u_j^d,\u_{j'}^d)\ge
\delta_\epsilon>0$ for any pair of critical points $\u_j^d,
\u_{j'}^d$ of $\A^\epsilon_d$.

Let $B_i$ be small isolating neighborhood of the critical points
$\u_i$ of $\A$. For $\epsilon>0$ sufficiently small all critical
points of $\A^\epsilon$ are contained in the neighborhoods $B_i$.
We can group together the critical points $\u_j^d$ of
$\A_d^\epsilon$ in associated neighborhoods $B_i^d$ in
$\overline{\DConf_d^n}$. Since the critical points $\u_i$ form a
Morse decomposition for $\A$ we obtain  the following Morse
inequalities from \cite[\S 7]{GVV}
\begin{equation}
\label{eq_Morse}
\sum_i P_\tau{(B_i^d)} = P_\tau(\HH) +(1+\tau)Q_\tau ,
\end{equation}
where $Q_\tau$ has nonnegative coefficients. Due to the uniform
separation of the critical points as $d\to \infty$
\rmref{eq_Morse} also holds in the limit for the functional
$\A^\epsilon$, i.e. $\sum_i P_\tau{(B_i)} = P_\tau(\HH)
+(1+\tau)Q_\tau$.

Lemma~\ref{CritGroup} below shows that for a given braid class,
each critical point $\u_i$ of $\A$ has Poincar\'e polynomial of
the form $P_\tau(\u_i)=A_i\tau^{p_i}$. By Lemma~\ref{MorsePert}
following, we can find Morse approximations $\A^\epsilon$ whose
Poincar\'e polynomial is exactly the same, i.e. $P_\tau(B_i) =
A_i\tau^{p_i}$. Substituting the latter into then Morse
inequalities for $A^\epsilon$ the proves that the number of
neighborhoods $B_i$ is bounded from below by the number of
monomials in $P_\tau(\HH)$ --- i.e. $|P_\tau(\HH)|$.
\qed

\begin{lemma}
\label{CritGroup} Given $\u$ an isolated critical point of $\A$,
the Poincar\'e polynomial is of the form $P_\tau(\u)=A\tau^p$ for
some $A\in \nat$ and $p\ge 0$.
\end{lemma}
\begin{proof}
In the case of a braid class with a single free strand, the
conclusion follows from a result of Dancer \cite{Dancer}: since
$\A$ is a first order Lagrangian of a scalar variable, a
degenerate critical strand has nullity at most two.

In the case of braids with multiple free strands, the proof
becomes somewhat more delicate. By considering the appropriate
covering we obtain an uncoupled system of equations for the
components of the braid $\u$. The critical groups of the braid
class are precisely the tensor product of the critical groups of
the individual components (see Theorem 5.5 of \cite{Chang}). Thus,
the Poincar\'e polynomials multiply, and the result follows from
the single-strand case.
\end{proof}

\begin{figure}[hbt]
\begin{center}
\psfragscanon
\psfrag{+}[][]{$+$}
\psfrag{-}[][]{$-$}
\includegraphics[width=5.25in]{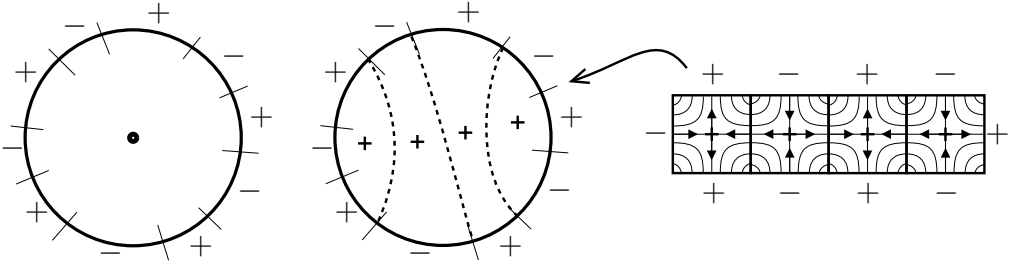}
\caption{[left] An isolating neighborhood of the critical point of $\A|_W$;
[right] Replace it with a chain of nondegenerate saddles.}
\label{fig_Isolating}
\end{center}
\end{figure}

\begin{lemma}
\label{MorsePert} Given $\A$ having finitely many critical points
$\u_i\in B_i$ with $P_\tau(\u_i)=A_i\tau^{p_i}$, there exists a
$C^2$-small perturbation of $\A$ with support in $\cup_iB_i$ to a
Morse functional $\A^\epsilon$ having exactly $A_i$ critical
points in $B_i$, each with Morse index $p_i$.
\end{lemma}
\begin{proof}
We consider each degenerate critical point separately. For each
degenerate critical point, the data in its critical groups comes
from a 2-dimensional `center' set $W$ given by the Gromoll-Meyer
version of the Morse Lemma \cite{GromMey}: all the
non-hyperbolicity of $d\A$ is manifested on $W$.

Consider $\A|_W:\real^2\to\real$ with coordinates chosen so that
there is a degenerate critical point at the origin having
$P_\tau=A_i\tau$. The statement of the lemma now becomes the claim
that there exists a perturbation of $\A|_W$ to a function on
$\real^2$ which has $A_i$ critical points of Morse index one. This
follows from choosing a small disc $D$ at the origin which is an
isolating neighborhood for $\nabla\A$. (This is possible via a
result of \cite{WilsonYorke}.) This implies that $\nabla\A$ is
transverse in/out of $\partial D$ on an alternating sequence of
$2A_i+2$ arcs as in Fig.~\ref{fig_Isolating}[left].

One may then set up analytic coordinates on $D$ and write out an
explicit Morse function with $A_i$ saddle points. A less explicit
method is to note that a linear chain of $A_i$ saddles --- as in
Fig.~\ref{fig_Isolating}[right] --- possesses an isolating
neighborhood whose boundary is combinatorially equivalent to that
of the disc $D$: for $D$ small, mapping this chain of saddles to
$D$ yields the appropriate perturbation of $\A$.
\end{proof}

\section{Proofs: Forcing periodic solutions}
\label{sec_proof2}

In this section, we provide details of the forcing arguments in
the case of non-stationary solutions. The technique is
philosophically the same as for stationary solutions: discretize,
apply the Morse-theoretic results of \cite{GVV}, then prove
convergence to solutions of \rmref{paracont}. However, the
requisite estimates are  more involved in the time-periodic case.
For this reason, we present the proofs for the normalized
equation,
\begin{equation}
\label{eq_normal}
    u_t = u_{xx}+g(x,u,u_x) ,
\end{equation}
noting that the general case of \rmref{paracont} is valid, though
messier.

Appendix~\ref{sec_appb} details a regularity result for
non-stationary solutions to \rmref{eq_normal}.

\subsection{Discretization and convergence}

We begin by truncating the system. Consider the equation
\begin{equation}
\label{eq_trunc}
    u_t=u_{xx}+g_K(x,u,u_x),
\end{equation}
where
\begin{equation*}
g_K(x,u,u_x) := \left\{
    \begin{array}{cl}
    g(x,u,u_x) & {\mbox{ for }} \abs{u}+\abs{u_x}\leq K \\
    \inf_{|u|+|u_x|\ge K}|g(x,u,u_x)| & {\mbox{ for }}
    \abs{u}+\abs{u_x}\geq K
    \end{array}\right. .
\end{equation*}
Consequently,
\[
    \abs{g_K(x,u,u_x)}\leq \abs{g(x,u,u_x)},
\]
for all $x\in S^1$, $u, u_x\in \real$. Thanks to this, the
estimates from Appendix~\ref{sec_appb} hold with the same
constants: any complete uniformly bounded solution $u^K(t,x)$ to
\rmref{eq_trunc} satisfies
\begin{equation}
\label{eq_truncsmooth}
    \abs{u_x^K}+\abs{u_{xx}^K}+\abs{u_{xxx}^K}+\abs{u_t^K}
    \leq C(\ell,
        \Vert u^K\Vert_{L^\infty}) ,
\end{equation}
with $C$ independent of the truncation domain $K$. By choosing $K$
appropriately, solutions of \rmref{eq_trunc} are also solutions of
\rmref{eq_normal}. Indeed, if $u^K(t,x)$ is a solution of
\rmref{eq_trunc} with $|u^K(t,x)|\le C_1$, then by
\rmref{eq_truncsmooth}, $|u^K_x(t,x)| \le C_2(\ell, C_1)$. If we
choose $K\ge \max{(C_1,C_2)}$, then solutions $u^K$ of
\rmref{eq_trunc}, with $|u^K(t,x)|\le C_1$, are also solutions of
\rmref{eq_normal}.

For convenience of notation we now omit the superscript $K$. We
discretize \rmref{eq_trunc} as follows: Let $u_i(t) = u(t,i/d)$
and
\begin{equation}
\label{eq_gkdisc}
    u_i' = d^2(u_{i+1}-2u_i+u_{i-1})
        + g_K\left(\frac{i}{d},u_i,d(u_{i+1}-u_i)\right)
        + \EE_i^d(u_{i-1},u_i,u_{i+1}),
\end{equation}
where $u_i'$ denotes $\frac{d}{dt}u(t,i/d)$. As before,
$|\EE_i^d|\le\abs{\epsilon_i(d)}\leq C/d$. The perturbations
$\EE_i^d$ are chosen such that the given stationary solutions of
\rmref{eq_normal} are also discretized solutions of
\rmref{eq_gkdisc}.

Let $\{u_i^d(t)\}$ be a sequence of solutions to \rmref{eq_gkdisc}
with $\abs{u_i^d(t)}\leq C_1$ for all $i$ and $d$. We will show
that one can pass to the limit as $d\to\infty$ and obtain a
complete solution to   \rmref{eq_normal}. The following lemma is
proved in a manner analogous to that of Lemma \ref{estI} of
\S\ref{sec_proof}.
\begin{lemma}
\[
    \int_J\sum_i\frac{1}{d}\abs{\d u_i}^2\,dt \leq C,
\]
where $J$ denotes the time interval $[T,T+1]$, and $C$ is
independent of $K$.
\end{lemma}
\begin{proof}
If we multiply \rmref{eq_gkdisc} by $u_i$ and then sum over
$i=0,..,d$ and integrate over $t\in [T,T+1]$ we obtain the desired
estimate as in Appendix \ref{sec_appb}. This uses the growth of
$g$ in $u_x$ given by Hypothesis {\bf (f2)}.
\end{proof}

Fix $K\ge\max{(C_1,C_2)}$, with $C_1$ and $C_2$ as above, and let
$f_i=g_K+\EE_i^d$. Then
\[
    \int_J\sum_i\frac{1}{d}\abs{f_i}^2\,dt\leq C .
\]
Write each solution $u_i(t)$ as a sum of terms $u_i=u_i^h +
u_i^p$, where
\begin{eqnarray*}
    \frac{d}{dt} u_i^h - \dd u_i^h = 0 , && u_i^h(T) = u_i(T) , \\
    \frac{d}{dt} u_i^p - \dd u_i^p = f_i , && u_i^p(T) = 0 .
\end{eqnarray*}
Then, for the homogeneous solutions $u_i^h$, one estimates
\[
    \int_{J'}\frac{1}{d}\sum_i\abs{\dd u_i^h}^2
    \leq
    {C\over d}\sum_i\abs{u_i(T)}^2\leq C,\quad J'=[T+\delta,T+1].
\]
This leads to the following estimate
$$
\int_{J'}\frac{1}{d}\sum_i\abs{(u_i^h)'}^2 dt
    + \int_{J'}\frac{1}{d}\sum_i\abs{\dd u_i^h}^2 dt
    \leq C .
$$
For the particular solution $u_i^p$, we have
$$
    \frac{1}{d}\sum_i\abs{f_i}^2=
    \frac{1}{d}\sum_i\abs{(u_i^p)'}^2
    -\frac{2}{d}\sum_i(u_i^p)'\dd u_i
    +\frac{1}{d}\sum_i\abs{\dd u_i^p}^2.
$$
For the middle term on the right hand side we have the identity
$
    -\frac{2}{d}\sum_i(u_i^p)'\dd u_i =
    \frac{d}{dt}\sum_i\frac{1}{d}\abs{\d u_i^p}^2
$. Upon integration over $J=[T,T+1]$ we obtain
$$
    \int_J\frac{d}{dt}\sum_i\frac{1}{d}\abs{\d u_i^p}^2\,dt=
    \frac{1}{d}\sum_i\abs{\left.\d u_i^p\right\vert_{T}^{T+1}}^2 =
    \frac{1}{d}\sum_i\abs{\d u_i^p(T+1)}^2 \ge 0.
$$
Combining these, we obtain
\[
    \int_J\frac{1}{d}\sum_i\abs{(u_i^p)'}^2 dt
    + \int_J\frac{1}{d}\sum_i\abs{\dd u_i^p}^2 dt
    \leq \int_J\frac{1}{d}\sum_i\abs{f_i}^2 dt
    \leq C.
\]
Combining the latter with the similar estimate for $u_i^h$ gives
the following estimate for the sum $u_i = u_i^p + u_i^h$:
\[
    \int_{J'}\frac{1}{d}\sum_i\abs{(u_i)'}^2 dt
    + \int_{J'}\frac{1}{d}\sum_i\abs{\dd u_i}^2 dt
    \leq C.
\]
Introduce the spline interpolation
\begin{eqnarray*}
\sp(u_i) &=& d\Delta^2 u_{i+1} (x-i/d)^3 - \Delta^2 u_{i+1} (x-i/d)^2\\
               &+& \Delta u_i (x-i/d) + u_i.
 \end{eqnarray*}
Now set $\widetilde U^d =\sp(u_i)$, and $U^d = \pl(u_i) = \Delta
u_i (x-i/d) + u_i$. Then,
\begin{eqnarray*}
\int_{J'}\int_{S^1} |\widetilde U^d -U^d|^2 dx dt
&\le& {C\over d^4}\to 0,\quad {\rm as}~~d\to \infty,\\
\int_{J'}\int_{S^1} |\widetilde U^d_x -U^d_x|^2 dx dt
&\le& {C\over d^2}\to 0,\quad {\rm as}~~d\to \infty,\\
\int_{S^1}|\widetilde U^d_{xx}|^2 dx
&\le & C \sum_i {1\over d}|\Delta^2 u_i|^2,\\
\int_{S^1}|\widetilde U^d_{t}|^2 dx
&\le & C \sum_i {1\over d}| u_i'|^2.
\end{eqnarray*}
From the latter two inequalities we derive that
$$
\widetilde U^d \in H^{1,2}(J';L^2(S^1))
\cap L^2(J';H^{2,2}(S^1)) \subset C(J';H^{1,2}(S^1)),
$$
which implies that
$
    \sum_i\frac{1}{d}\abs{\d u_i(t)}^2 \leq C
    \quad   \forall t\in\real .
$
Moreover,
\begin{eqnarray*}
\widetilde U^d, U^d &\rightarrow&
    u, \quad {\rm in}~~L^2(J';H^{1,2}(S^1)),\\
\widetilde U^d_t, U^d_t &\rightharpoonup &
    u_t \quad {\rm in}~~L^2(J';L^2(S^1)).
\end{eqnarray*}
From these embeddings one easily deduces that
$$
g_K(x,U^d,U_x^d) \longrightarrow g_K(x,u,u_x),
    \quad {\rm in}~~L^2(J';L^2(S^1)).
$$
Choose smooth test functions of the form $\phi(t,x) = \sum_{k=1}^N
\alpha_k(t) w_k(x)$, where $\{w_k\}$ is an orthonormal basis for
$H^{1,2}(S^1)$. Set $\phi_i(t) = \phi(t,i/d)$, and $\Phi^d =
\pl(\phi_i)$, then
$$
\int_{J'} \sum_i {1\over d}\Bigl( g_K(i/d,u_i,\Delta u_i)
+\EE_i^d\Bigr) \phi_i dt\longrightarrow
\int_{J'}\int_{S^1} g_K(x,u,u_x) \phi~ dx dt.
$$
Because of the PL approximation the
following  integrals become sums over the anchor points:
\begin{eqnarray*}
    \int_{S^1} U^d_x\Phi_x dx
&=&
    \sum_i\frac{1}{d}\d u_i\d\phi_i
    = -\sum_i\frac{1}{d}\dd u_i\phi_i ,
\\
    \int_{S^1} U^d_t\Phi\,dx
&=&
    \sum\frac{1}{d}u_i' \phi_i
    + \frac{1}{3d}\sum_i\frac{1}{d}(u'_{i+1}-u'_i)\d\phi_i ,
\\
    \int_{S^1} f\Phi\,dx
&=&
    \sum_i\frac{1}{d}f_i\phi_i
    + \frac{1}{2d}\sum_i\frac{1}{d}f_i\d\phi_i .
\end{eqnarray*}
The final terms of the last two equations admit the
following bounds:
\vfill
\begin{eqnarray*}
    \abs{
    \frac{1}{3d}\sum_i\frac{1}{d}(u'_{i+1}-u'_i)\d\phi_i
    }
&
    \leq
&
    \frac{2}{3d}
    \left(\int_J\sum_i\frac{1}{d}\abs{f_i}^2dx\right)^{\frac{1}{2}}
    \left(\int_J\sum_i\frac{1}{d}\abs{\d\phi_i}^2dx\right)^{\frac{1}{2}}
\\
&
    \leq
&
    \frac{C}{d} \to 0 ,
\\
    \frac{1}{2d}\sum_i\frac{1}{d}f_i\d\phi_i
&
    \leq
&
    \frac{1}{2d}
    \left(\int_J\sum_i\frac{1}{d}\abs{u'_i}^2dx\right)^{\frac{1}{2}}
    \left(\int_J\sum_i\frac{1}{d}\abs{\d\phi_i}^2dx\right)^{\frac{1}{2}}
\\
&
    \leq
&
     \frac{C}{d} \to 0 .
\end{eqnarray*}
Weak convergence implies that as $d\to\infty$,
\[
    \int_J\int_{S^1}\Bigl[ U^{d}\Phi
+    U^{d}_x\Phi_x \Bigr]\,dx\,dt
\longrightarrow
    \int_J\int_{S^1} \Bigl[ u\phi
+    u_x\phi_x \Bigr]\,dx\,dt,
\]
$$
\int_J\int_{S^1} U_t^{d}\Phi \,dx\,dt
\longrightarrow
    \int_J\int_{S^1}  u_t\phi \,dx\,dt,
$$
where $u(t,x)$ is the weak limit of $U^d(t,x)$.
Hence, $u$ is a weak solution to \rmref{eq_trunc} for all smooth
test function  $\phi$ defined above. These functions form a dense
subset in $ H^{1,2}(J'\times S^1)$, and therefore, since $u_i$
satisfies \rmref{eq_gkdisc},
\[
    \int_{S^1} u_t\phi\,dx
+   \int_{S^1} u_x\phi_x\,dx
=   \int_{S^1} g_K(x,u,u_x)\phi\, dx ,\quad
        \forall~\phi \in H^{1,2}(S^1).
\]
Standard regularity theory arguments then yield strong solutions
to \rmref{eq_trunc}. The using the $L^\infty$-bounds on $u$ we
also conclude that $u$ is a weak solution to \rmref{eq_normal}.
Using standard regularity techniques one can show that the
convergence is in $H^{1,2}(J'\times S^1)$. This completes the
proof of the following theorem:
\begin{theorem}
\label{thm_conv-per} For any sequence of bounded solutions
$\{u_i^d(t)\}$ of \rmref{eq_gkdisc} with $|u_i^d(t)| \le C$, for
all $t$ and $i$, $\pl(u_i^d)$ converges, in $H^{1,2}(J\times
S^1)$, to a (strong) solution $u$ of \rmref{eq_trunc}. Moreover if
$K$ is chosen large enough then $u$ is a (strong) solution of
\rmref{eq_normal}.
\end{theorem}

\subsection{Proof of Theorem~\ref{thm_Main2}}

Let $\{\u\rel\v\}$ be a braid class that does not permit
stationary solutions for \rmref{eq_normal}. For $d$ large enough
the same holds for \rmref{eq_gkdisc}; otherwise, the results in \S
\ref{sec_proof} would yield stationary solutions of
\rmref{eq_normal}, a contradiction. If $\{\u\rel\v\}$ is bounded
and proper with $\HH(\u\rel\v)\not = 0$, then for each $d$ large
enough there exists a periodic  solution $\uu^d$ with strands
$u_i^{\alpha,d}(t)$ via \cite[Thm. 2]{GVV}. By Theorem
\ref{thm_conv-per} this sequence yields a solution $u(t,x)$ of
\rmref{eq_normal}.

It remains to be shown that $u(t,x)$ is periodic in $t$. This
follows from the celebrated Poincar\'e-Bendixson Theorem for
scalar parabolic equations due to Fiedler and Mallet-Paret
\cite{FMP}, which states that a bounded solution $u(t,x)$ has
forward limit set either a stationary point or a time-periodic
orbit. By assumption $\{\u\rel\v\}$ contains no stationary points
which leaves the second option;  a periodic solution. This also
proves then that $\{\u\rel\v\}$ contains a periodic solution of
the desired braid class.
\qed

We remark that the proof above is for braid classes $\{\u\rel\v\}$
for which $\u$ has a single component. For $\u$ with multiple
components, a nonvanishing index implies that each component of
$\u$ is either stationary or periodic; however, unless the periods
are rationally related, the entire braid class will be merely
quasi-periodic as opposed to periodic.

\section{Concluding remarks}\label{sec_conc}

{\it Boundary conditions.} We have employed periodic boundary
conditions for convenience and as a means to allow for
time-periodic orbits. Nothing prevents us from using other
boundary conditions, although the resulting dynamics is often
gradient-like. Neumann, Dirichlet, or (nonlinear) combinations of
the two are imposed by choosing closed subsets
$B_0\subset\{(0,u,u_x)\}$ and $B_1\subset\{(1,u,u_x)\}$ and
requiring the braid endpoints to remain in these subspaces. As the
topology of the configuration spaces of braids may change, so may
the resulting invariants. Since the comparison principle still
holds, our topological methods remain valid, though the invariants
themselves may change.

{\it Coercivity and unbounded classes.} Theorem \ref{thm_infinite}
deals with dissipative systems. The opposite of dissipative is the
{\em coercive} condition:
$$
uf(x,u,0,0) \to \infty,\quad {\rm as}~~|u| \to \infty,
$$
for all $x \in S^1$. For either of these cases, the restriction to
bounded braid classes may be relaxed. For dissipative systems, any braid
class becomes bounded by adding two unlinked strands as per
Appendix~\ref{appendixC}. In order to deal with coercive systems one needs to
include the behavior of the system at infinity.
We propose that a compactification of the unbounded braid classes
yields an index with the same properties as that for bounded classes.

{\it Improper braids.} A braid class is improper if components of
the braid can be collapsed. Our results on $t$-periodic solutions
in \S\ref{sec_exper} dealt with improper braids in an ad hoc
manner by  `blowing up' the collapsible strands via adding
additional strands to the skeleton.

A different approach would be to blow up the vector field in the
traditional manner via homogeneous coordinates, working in the
setting of finite-dimensional PRRs. Stabilization then allows one
to define the invariant in the continuous limit. This type of
blow-up procedure is very general and should be applicable to a
wide variety of systems.

{\it Periodic skeleta.} The forcing theory we have developed uses
stationary solutions for the skeleton. We believe that all of the
results hold for skeleta composed of time-periodic orbits.

{\it p-Laplacians and degenerate parabolic equations.} The fully
nonlinear parabolic equations studied in this paper are restricted
by the `uniform parabolicity' hypothesis given by {\bf (f1)}. We
choose to restrict ourselves to uniform parabolic equations in
order to keep technicalities to a minimum. However, the theory
should also apply to degenerate parabolic equations of various
kinds. One weakening of Hypothesis {\bf (f1)} would read

\indent\indent$\quad 0<\partial_w f(x,u,v,w),\quad$ ~~~~~~for all
    ~~$ w\not=0$, and $(x,u,v) \in S^1\times\real^2$.

Good examples of degenerate equations are the 1-dimensional porous
medium equation $u_t = (u^p u_x)_x + g(x,u,u_x)$, or the
p-Laplacian equation $u_t = (|u_x|^{p-1}u_x)_x + g(x,u,u_x)$.
Solutions of these equations have less regularity than
\rmref{paracont}, which complicates the approach used in
\S\ref{sec_proof}. In that case, one can use the weak solution
approach as carried out in the periodic case. The key point is to
find the appropriate estimates in $u_x$.

{\it Scalar hyperbolic conservation laws.} Conservation laws of
the form
\begin{equation}\label{eq_hyp}
u_t = f(x,u,u_x),
\end{equation}
where $f$ is monotonically increasing in $u_x$, discretize to
one-sided parabolic systems of the form $u_i' =
\R_i(u_i,u_{i+1})$, cf. \cite{MPSm}. Our theory remains valid for
discretized systems of this form; if we establish the appropriate
a priori $L^\infty$-estimates a braid-forcing theory for
\rmref{eq_hyp} can be derived.

\appendix

The following estimates, though necessary, are antithetical to our
philosophy: the entire forcing theory for \rmref{paracont} is
topological in nature.

\section{Estimates: stationary}
\label{sec_appa}

A stationary solution of \rmref{paracont} is some $u \in
C^2(\re/\ell \zed)$  satisfying $f(x,u,u_x,u_{xx})=0$. Hypotheses
{\bf(f1)}-{\bf(f2)} permit the following regularity statement.
\begin{lemma}\label{lem_A}
Let $u\in C^2(\re/\ell \zed)$ be a stationary solution of
\rmref{paracont} with $f$ satisfying {\bf(f1)}-{\bf(f2)}. There
exists a constant $C=C(\ell, \norm{u}_{L^\infty})$ depending only
on the sup-norm of $u$, such that
\begin{equation}
\label{eq_normsa}
    \abs{u_x}+\abs{u_{xx}}+\abs{u_{xxx}} \le C.
\end{equation}
\end{lemma}
\begin{proof}
Using {\bf(f1)} we obtain the following estimate for $f$;
\begin{equation}\label{paracontA}
a_-(u_{xx}) u_{xx} + f(x,u,u_x,0) \le f(x,u,u_x,u_{xx}) \le
a_+(u_{xx}) u_{xx} + f(x,u,u_x,0),
\end{equation}
where $a_-$ and $a_+$ are defined in  \S\ref{sec_proof}. Multiply
\rmref{paracontA} by $u$. Integrating over $S^1:=\real/\ell \zed$,
using Hypothesis {\bf(f2)} and the fact that ${1\over a_\pm} \le
\lambda^{-1}$ yields
\begin{eqnarray*}
\int_S u_x^2 dx &\le& \int_S \lambda^{-1} |u|\cdot |f(x,u,u_x,0)| dx\\
&\le&
C\norm{u}_{L^\infty} \int_S|f(x,u,u_x,0)|dx\\
&\le& C
    \left(1+
    \int_S|u_x|^\gamma dx
    \right) .
\end{eqnarray*}
Since $\gamma<2$, it follows that $\int_S |u_x|^2 dx \le C$. Next
we deduce from \rmref{paracontA} that $|u_{xx}| \le \lambda^{-1}
|f(x,u,u_x,0)|$. Again by using Hypothesis {\bf(f2)} we obtain
\begin{eqnarray*}
\int_S |u_{xx}|^{2\over \gamma}
&=&
C\int_S|f(x,u,u_x,0)|^{2\over \gamma} dx\\
&\le&
\int_S\Bigl| C +
    C|u_x|^\gamma\Bigr|^{2\over\gamma} dx\\
&\le& C
    \left(1+
    \int_S|u_x|^2 dx
    \right)
    \le C.
\end{eqnarray*}
The latter implies that $\norm{u}_{W^{2,{2\over\gamma}}} \le C$.
From the Sobolev embeddings for $W^{2,{2\over \gamma}}(S)$ we
derive
$$
\norm{u}_{C^{1,\alpha}(S)} \le C \norm{u}_{W^{2,{2\over\gamma}}}
\le C,
$$
with $0<\alpha<1-{\gamma\over 2} <1$. In particular $\Vert
u_x\Vert_{L^\infty} \le C$. Again by using  the pointwise bound
$|u_{xx}| \le \lambda^{-1} |f(x,u,u_x,0)|$ we obtain
\begin{eqnarray*}
\sup_{x} |u_{xx}|
&\le &
C\Vert f(x,u,u_x,0)\Vert_{L^\infty}\\
&\le&
C
    + C\Vert u_x \Vert^\gamma_{L^\infty}
\le
C,
\end{eqnarray*}
which implies that $\norm{u_{xx}} \le C$. By differentiating the
equation and using the fact that $f\in C^1$ to estimate $u_{xxx}$,
we obtain
\[
\partial_x f + \partial_u f\cdot u_x + \partial_{u_x}  f \cdot u_{xx}
+ \partial_{u_{xx}} f \cdot u_{xxx} =0 .
\]
 For $u_{xxx}$ this yields
\begin{eqnarray*}
|u_{xxx}| &\le& {1\over \partial_{u_{xx}} f} \Bigl\{ |\partial_x
f| + |\partial_{u_x} f||u_x| + |\partial_{u_{xx}} f||u_{xx}|
\Bigr\} \,\le\, C,
\end{eqnarray*}
since all derivatives of $f$ can be bounded in terms of
$\norm{u}_{L^\infty}$. This completes the proof.
\end{proof}

\section{Estimates: non-stationary}
\label{sec_appb}

We repeat the regularity arguments for non-stationary solutions to
\rmref{paracontII}. As the estimates are similar in spirit as
those of Appendix~A, we omit the more unseemly steps.
\begin{lemma}\label{lem_B}
Let $u\in C^1(\real;C^2(\real/\ell \zed))$ be a complete bounded
solution with $g$ satisfying {\bf(f1)}-{\bf(f2)}. There exists a
constant $C=C(\ell, \norm{u}_{L^\infty})$ depending only on the
sup-norm of $u(t,x)$, such that
\begin{equation}
\label{eq_normsb}
    \abs{u_x}+\abs{u_{xx}}+\abs{u_{xxx}} +\abs{u_t}\le C.
\end{equation}
\end{lemma}
\begin{proof}
As before, let $S^1:=\real/\ell \zed$. Denote by $J$ the time
interval $J:=[T,T+1]$. Multiplying \rmref{paracontII} by $u$ and
integrating by parts yields
\[
    \int_J\int_{S^1}u_tu\,dx\,dt
    =
    -\int_J\int_{S^1}u_x^2\,dx\,dt
    +\int_J\int_{S^1}g(x,u,u_x)u\,dx\,dt .
\]
Using hypothesis {\bf(f2)} we derive
\[
    \int_J\int_{S^1} u_x^2\,dx\,dt
    \leq
    -\frac{1}{2}\left.\int_{S^1}u^2\,dx\right\vert_T^{T+1}
    +C+C\int_J\int_{S^1}\abs{u_x}^\gamma\,dx\,dt.
\]
Hence, since $\gamma<2$, $\int_J\int_{S^1}|u_x|^2\,dx\,dt \le C$.

We proceed with the more technical estimates. Given the solution
$u(t,x)$,
\[
    u_t-u_{xx} = g(x,u(t,x),u_x(t,x))
        \in L^{\frac{2}{\gamma}}(J;L^{\frac{2}{\gamma}}(S^1)) ,
\]
since $|g|^{2/\gamma} \le C +C|u_x|^2$. As such, $L^p$ regularity
theory implies (see, e.g., \cite{Brezis})
\begin{eqnarray*}
    \norm{u_t}_{L^{\frac{2}{\gamma}}(J';L^{\frac{2}{\gamma}}(S^1))}
        &\leq&
    C(\delta)\norm{f}_{{L^\frac{2}{\gamma}}(J;L^{\frac{2}{\gamma}})},
        \\
    \norm{u_{xx}}_{L^{\frac{2}{\gamma}}(J';L^{\frac{2}{\gamma}}(S^1))}
        &\leq&
    C(\delta)\norm{f}_{{L^\frac{2}{\gamma}} (J;L^{\frac{2}{\gamma}})},
\end{eqnarray*}
where $J':=[T+\delta,T]\subset J$ for some $0<\delta\ll 1$. In
particular,
\[
    u\in
    L^{\frac{2}{\gamma}}(J';H^{2,\frac{2}{\gamma}}(S^1))
    \cap
    L^\infty(\real;L^\infty(S^1)) .
\]
Bootstrapping proceeds in a standard fashion using a parabolic
version of the Gagliardo-Nirenberg interpolation inequalities.
Given any function $u\in L^p(J',H^{2,p}(S^1))\cap
L^\infty(J',L^\infty(S^1))$, then
\[
    \norm{u}_{L^{2p}(J';H^{1,{2p}}(S^1))}
    \leq
    C\norm{u}^{\frac{1}{2}}_{L^p(J',H^{2,p})}
    \cdot
    \norm{u}^{\frac{1}{2}}_{L^\infty(J',L^\infty)} .
\]
Therefore, we have $u\in
L^{\frac{4}{\gamma}}(J',H^{1,\frac{4}{\gamma}}(S^1))$ and, hence,
$g\in L^{\frac{4}{\gamma^2}}(J';L^{\frac{4}{\gamma^2}}(S^1))$.

We repeat the procedure $k$ times, each time restricting the time
domain $[T+k\delta,T+1]$. Choose $k>0$ sufficiently large so that
$(2/\gamma)^k>2$ and choose $\delta$ sufficiently small so that
$[T+k\delta,T+1]$ contains $J'':=[T+\frac{1}{2},T+1]$. Then we
have
\[
    u\in
    H^{1,2}(J'';L^2(S^1))
    \cap
    L^2(J'';H^{2,2}(S^1)) .
\]
By Sobolev embedding, we get $u\in C(J'';H^{1,2}(S^1))$. Repeating
the entire procedure yields $u\in
C^\alpha(J'';C^{1,\alpha}(S^1))$. This bound is now independent of
$T$, and one translates to obtain $u\in
C^\alpha(\real;C^{1,\alpha}(S^1))$. The additional smoothness now
follows directly from the fact that $u$ solves \rmref{paracontII}.

The $C^3$-estimate is obtained as in the stationary case by
differentiating the equation and using the $C^{1,2}$-estimates
obtained above.
\end{proof}

\section{Discrete enclosure}
\label{appendixC} Using a discrete version of enclosure between
sub/super solutions and a nontrivial braid diagram, we obtain the
following existence result.
\begin{lemma}\label{help}
Let $f$ satisfy Hypotheses {\bf(f1)}-{\bf(f2)} and let $\v$ be a
non-trivially braided stationary braid for \rmref{paracont}.
Assume that there exists a   $u^*$ such that $v^\alpha(x) < u^*$
for all $\alpha$ and $f(x,u^*,0,0)<0$. Then, there exists a
1-periodic solution $u$ with
$$
\max_{\alpha}v^\alpha(x) < u(x) < u^*,
$$
for all $x\in S^1$.
\end{lemma}
It is clear that the result holds for case of a  $u^*$ such that
$u^* < v^\alpha(x)$ for all $\alpha$ and $f(x,u^*,0,0)>0$. In that
case one finds a solution $u$ satisfying
$$
u^* < u(x) < \min_{\alpha}v^\alpha(x)
$$
for all $x\in S^1$.
\proof As in \S\ref{sec_proof} we discretize
\rmref{paracont} in $x$. For $u^*$ this implies that
$f(i/d,u^*,0,0)<0$. As for the braid $\v$ we use Lemma
\ref{approx} to find $\EE_i^d$ and the recurrence relation
$\R_i^d(u_{i-1},u_i,u_{i+1}) :=
    f(i/d,u_i,\d u_i, \dd u_i) + \EE_i^d(u_{i-1},u_i,u_{i+1})$.
By construction the discretized skeleton $\disc_d\v$ is stationary
for $\R$.

Define the region
$$
D = \{ \{u_i\}_{i=0}^d~|~\max_{\alpha} v_i^\alpha \le u_i \le
u^*,~~u_0=u_d\}.
$$
If the discretization is chosen fine enough then the discretized
braid is non-trivial. As a consequence $u_i$ cannot collapse onto
$\disc_d\v$ and if $u_i = v^\alpha_i$ for some $i$ and some
$\alpha$, then $\R_i^d(u_{i-1},u_i,u_{i+1})>0$. By the definition
of $u^*$ it follows that if $u_i =u^*$ for some $i$, then
$\R_i^d(u_{i-1},u_i,u_{i+1}) \le f(i/d,u_i,0,0)<0$ (parabolicity).
The region $D$ is therefore an attracting isolating (compact) set
for \rmref{eq_PRR}. Thus for each large enough $d$ we find a
discrete solution $\{u_i^d\}_{i=0}^d$. Since $\{u_i^d\}_{i=0}^d$
is a priori bounded we derive from the limiting procedure in
\S\ref{sec_proof} that this yields a stationary solution $u(x)$
for \rmref{paracont}, satisfying the desired inequality. \qed


\begin{thebibliography}{9}
{\footnotesize

\bibitem{Ang1} S.B. Angenent, \emph{ The zero set of a solution of a
parabolic equation}, J. Reine Ang. Math. {\bf 390}, 1988, 79-96.

\bibitem{Ang2} S.B. Angenent, \emph{Curve Shortening and the topology
of closed geodesics on surfaces}, preprint 2000.

\bibitem{Ang3} S.B. Angenent, \emph{The periodic orbits of an area
preserving twist map}, Comm. Math. Phys. {\bf 115}, 1988, 353-374.

\bibitem{AngFied} S.B. Angenent and B. Fiedler, \emph{ The dynamics of
rotating waves in scalar reaction diffusion equations}, Trans. AMS
{\bf 307}(2) 1988, 545-568.

\bibitem{Bir} J. Birman, \emph{Braids, links and the mapping class group}, Ann.
Math. Stud. {\bf 82} 1975, Princeton Univ. Press.

\bibitem{Brezis} H. Brezis, \emph{Analyse Fonctionnelle}, Mason, 1983.

\bibitem{Chang} K.C. Chang, \emph{Infinite Dimensional Morse Theory and
Multiple Solution Problems}, Birkhauser 1991.

\bibitem{Conley} C. Conley, \emph{Isolated Invariant Sets and the Morse Index},
CBMS Reg. Conf. Ser. Math. {\bf 38} 1978, published by the AMS.

\bibitem{Dancer} N. Dancer, \emph{ Degenerate critical points, homotopy
indices and Morse inequalities}, J. Reine Ang. Math. {\bf 350} 1984, 1-22.

\bibitem{FMP} B. Fiedler and J. Mallet-Paret, {\emph A Poincar\'e-Bendixson
theorem for scalar reaction diffusion equations.}
Arch. Rational Mech. Anal.  {\bf 107}, 1989,  no. 4, 325--345.

\bibitem{FOl} G. Fusco and W. Oliva, \emph{ Jacobi matrices and
transversality}, Proc. R. Soc. Edinb. Sect. A Math.
{\bf 109} 1988, 231-243.

\bibitem{GVV} R.W. Ghrist, J.B. Van den Berg and R.C. Vandervorst,
\emph{ Morse theory on spaces of braids and Lagrangian dynamics},
Invent. Math. {\bf 152} 2003, 369-432.

\bibitem{GromMey} D. Gromoll and W. Meyer, {\emph On differentiable
functions with isolated critical points.}  Topology {\bf 8},  1969, 361--369.

\bibitem{Hale} J.K. Hale, {\emph Dynamics of a scalar parabolic equation.}
Canad. Appl. Math. Quart. {\bf 5}, 1997, no. 3, 209--305.

\bibitem{MPSm} J. Mallet-Paret and H.L. Smith, \emph{ The Poincar\'e-Bendixson
theorem for monotone cyclic feedback systems}, J. Dyn. Diff. Equations
{\bf 2}, 1990, 367-421.

\bibitem{Matano1} H. Matano, \emph{ Nonincrease of the lap-number of a
solution for a one-dimensional semi-linear parabolic equation},
J. Fac. Sci. Tokyo 1A {\bf 29} 1982, 645-673.

\bibitem{Nakashima1} K. Nakashima, \emph{Stable transition layers in a
balanced bistable equation}.  Diff. Integral Equations {\bf 13}, 2000,
no. 7-9, 1025--1038.

\bibitem{Nakashima2} K. Nakashima, \emph{Multi-layered stationary solutions
for a spatially inhomogeneous Allen-Cahn equation.}
J. Differential Equations  {\bf 191},  2003,  no. 1, 234--276.

\bibitem{Smillie} J. Smillie, \emph{ Competative and cooperative tridiagonal
systems of differential equations}, SIAM J. Math. Anal.
{\bf 15} 1984, 531-534.

\bibitem{Sturm} C. Sturm, \emph{ M\'emoire sur une classe d'\'equations
\`a diff\'erences partielles},
J. Math. Pure Appl. {\bf 1} 1836, 373-444.

\bibitem{WilsonYorke} W. Wilson and J. York, \emph{Lyapunov functions
and isolating blocks,}
J. Differential Equations {\bf 13}, 1973, 106--123.

\bibitem{Zel} T. I. Zelenyak, \emph{Stabilization of solutions of boundary value
problems for a second order parabolic equation with one space
variable}, Differential equation {\bf 4},
                    1968, 27-22.


}\end{thebibliography}
\end{document}